\documentclass[10pt]{article}

\usepackage{latexsym}
\usepackage{amsfonts}
\usepackage{amsmath}
\usepackage{mathrsfs}  
\usepackage{dsfont}    
\usepackage{bbold}     
\usepackage[english]{babel}
\usepackage{caption}
\usepackage{epsfig}
\usepackage{float}
\usepackage{isolatin1}
\usepackage{subfigure}

\newtheorem{theorem}{Theorem}

\newtheorem{lemma}{Lemma}

\newcommand{\zerarcounters}{\setcounter{equation}{0}}

\newcommand{\ZZZ}{\mathds{Z}}
\newcommand{\CCC}{\mathds{C}}
\newcommand{\NNN}{\mathds{N}}

\newcommand{\RRR}{\mathds{R}}
\newcommand{\TTT}{\mathds{T}}
\newcommand{\one}{\mathds{1}}

\newcommand{\BB}{{\mathcal B}}

\newcommand{\MM}{{\mathcal M}}

\newcommand{\PP}{{\mathcal P}}

\newcommand{\RR}{{\mathcal R}}
\newcommand{\SSS}{{\mathcal S}}

\newcommand{\VV}{{\mathcal V}}

\newcommand{\gotl}{{\mathfrak l}}
\newcommand{\gotm}{{\mathfrak m}}

\newcommand{\gots}{{\mathfrak s}}

\newcommand{\gotM}{{\mathfrak M}}

\newcommand{\esselle}{{\gots\gotl}}
\newcommand{\SL}{{\rm  SL}}
\newcommand{\tr}{{\rm tr}\,}
\newcommand{\Val}{{\rm Val}}

\newcommand{\Fullbox}{{\rule{2.0mm}{2.0mm}}}

\newcommand{\EP}{\hfill\Fullbox\vspace{0.2cm}}
\newcommand{\prova}{\noindent{\it Proof. }}

\newcommand{\eps}{\varepsilon}
\newcommand{\al}{\alpha}

\newcommand{\ka}{\kappa}
\newcommand{\g}{\gamma}
\newcommand{\om}{\omega}
\newcommand{\la}{\lambda}

\newcommand{\oo}{\boldsymbol{\omega}}
\newcommand{\nn}{\boldsymbol{\nu}}
\newcommand{\pps}{\boldsymbol{\psi}}
\newcommand{\vzero}{\boldsymbol{0}}

\def\ins#1#2#3{\vbox to0pt{\kern-#2 \hbox{\kern#1 #3}\vss}\nointerlineskip}

\begin{document}

\title{\bf Resummation of perturbation series and\\
reducibility for Bryuno skew-product flows}
\author
{\bf Guido Gentile
\vspace{2mm} \\ \small 
Dipartimento di Matematica, Universit\`a di Roma Tre, Roma,
I-00146, Italy.
\\ \small 
E-mail: gentile@mat.unirom3.it
}

\date{}

\maketitle

\begin{abstract}
We consider skew-product systems on $\TTT^{d} \times \SL(2,\RRR)$ for
Bryuno base flows close to constant coefficients, depending on a
parameter, in any dimension $d$, and we prove reducibility for a large
measure set of values of the parameter. The proof is based on a
resummation procedure of the formal power series for the conjugation,
and uses techniques of renormalisation group in quantum field theory.
\end{abstract}




\zerarcounters
\section{Introduction}
\label{sec:1}

Consider the linear differential equation
\begin{equation}
\dot x = \left( \la A + \eps f(\oo t) \right) x ,
\label{eq:1.1} \end{equation}
on $\SL(2,\RRR)$, where $\la\in[a,b]\subset\RRR$,
$\eps$ is a small real parameter, $\oo\in\RRR^{d}$ is a vector with
rationally independent components, and $A,f\in \esselle(2,\RRR)$,
with $A$ is a constant matrix and $f$ an analytic function periodic
in its arguments. We say that $f$ is quasi-periodic in time $t$.

Reducibility for (\ref{eq:1.1}) means the existence of a
quasi-periodic change of variables which takes the system into
a system with constant coefficients:
\begin{equation}
x=B(\oo t)y, \qquad \dot y = A_{0} y ,
\label{eq:1.2} \end{equation}
with $B \in \SL(2,\RRR)$ analytic and $A_{0} \in \esselle(2,\RRR)$
constant. In particular if the solution $y(t)$ is periodic then the
solution $x(t)$ is quasi-periodic, hence bounded for all times.

A special case of (\ref{eq:1.1}) is the one-dimensional
Schr\"odinger equation with a weak quasi-periodic potential,
or with arbitrary quasi-periodic potential for large energy.
By assuming a suitable non-resonance condition on the frequency vector
$\oo$, reducibility for $\eps$ small enough and for a large measure set
of values $\la$ in $[a,b]$ (for which quasi-periodic solutions exist)
was proved by Dinaburg and Sinai \cite{DS}, by using KAM techniques;
see also \cite{PF} for a review. Weaker non-resonance conditions
were shown to be possible by R\"ussmann \cite{R}, then used by
Moser and P\"oschel \cite{MP} to enlarge the set of values $\la$ for
which reducibility can be obtained. Reducibility almost everywhere
in $\la$ and for small $\eps$ has been obtained by Eliasson \cite{E2},
for $\oo$ a Diophantine vector.

A brief survey on the problem of reducibility for skew-product
systems can be found in \cite{E3,E4}. In particular results similar
to those by Eliasson, -- i.e. reducibility almost everywhere for
Diophantine frequency vectors, -- in the case of other Lie groups,
also not close to constant coefficients, have been obtained by
Krikorian\cite{K1,K2}. Very recently, Avila and Krikorian \cite{AK}
proved, by using renormalization techniques, that, if $\oo$ belongs
to a subset of full measure of the Diophantine vectors in $d=2$,
for all values of $\eps$ and almost everywhere in $\la$, quasi-periodic
Schr\"odinger cocycles are either reducible or non-uniformly hyperbolic.

R\"ussmann's non-resonance condition is weaker than the usual
Diophantine one, and is expressed in terms of a suitable
approximation function \cite{R,MP}. In $d=2$ it is equivalent to
Bryuno's condition. Bryuno vectors have been explicitly considered in
the case of skew products for the first time by Lopes Dias in \cite{L},
where in $d=2$ a normal form theorem (analogous to Lemma \ref{lem:22}
below) is proved with renormalisation group techniques. However,
the non-resonance condition with the eigenvalues of $\la A$
is still assumed to be of Diophantine type. Renormalisation group
techniques have been also used in \cite{KL} for any $d$,
in the case of Diophantine vectors.

In this paper we consider Bryuno vectors in any dimension, and,
for $\eps$ small enough, we prove reducibility on a large measure
set of values of $\la$. The estimates we find for the excluded set
are much better than those provided by standard KAM methods (cf. for
instance \cite{K1}). The techniques we use are those of
renormalisation group typical of quantum field theory,
based on a diagrammatic representation of the equation in terms of
trees, and are inspired to those used in \cite{G1,G2,GG2}.
Trees for skew-products were already introduced by Iserles
and N{\o}rsett \cite{IN,I}, but they used expansions in time,
hence not suited for the study of global properties,
such as reducibility and quasi-periodicity.

We can formulate our result as follows.

\begin{theorem}\label{thm:1}
Let $A\in\esselle(2,\RRR)$ be a constant matrix with imaginary
eigenvalues and $f\in\esselle(2,\RRR)$ an analytic quasi-periodic
function of time. Let $\oo\in\RRR^{d}$ be a Bryuno vector.
Then there exists $\eps_{0}>0$ and $\sigma>0$ such that
for all $|\eps|<\eps_{0}$ the set of values $\la\in [a,b]$
for which the system (\ref{eq:1.1}) is not reducible
is of Lebesgue measure less than ${\rm const.}|\eps|^{\sigma}$.

\end{theorem}

The proof will proceed through the following steps.
In Section \ref{sec:2}, we reduce the study of system (\ref{eq:1.1}) to
the study of a system of differential equations in $\CCC^{2}$,
that we call here the ``auxiliary system'', and we see that
the property for $x$ to have $\det x=1$ can be interpreted as
the existence of a suitable first integral for the new system.
Next, in Section \ref{sec:3} we look for a quasi-periodic solution
of the auxiliary system: first, we try for solutions in the form
of formal power series in $\eps$. However, in order
to define such a series, even order by order, we cannot fix $\la$.
Instead, we write $\la=\la_{0}+\mu$, with $\la_{0}$ in some
interval $\Lambda_{0}$, and we see that for fixed $\la_{0}$
there exists a formal power series for $\mu$ such that
the auxiliary system admits a formal power series solution.
Hence the formal series turn out to be well-defined order by order.
Moreover, they can be proved to be convergent. This result can be
obtained by performing a suitable resummation leading to simpler series,
that we shall call renormalised series (for a discussion of the method
within the standard KAM theory we refer to \cite{GG1,GBG}).
The renormalised series are introduced in Section \ref{sec:4}, and
in Section \ref{sec:5} are proved to converge to well-defined functions.
The latter are analytic in $\eps$ and solve the differential equation
of the auxiliary system, provided $\lambda_{0}$ is chosen in a
subset $\Lambda_{0}^{*}$ of $\Lambda_{0}$.
Finally in Section \ref{sec:6} we have to control that
the set of values $\la\in[a,b]$ for which the above procedure
can be followed coincide with $[a,b]$, up to a small measure set.

We conclude with two comments.

Given the system (\ref{eq:1.1}) one could also consider another problem:
fix $\la$ and study for which values of $\eps$ (small enough) the system
is reducible. This a natural question if, for instance, 
instead of the Schr\"odinger equation, one considers Hill's equation,
where there is no free parameter other than $\eps$ itself.
Under suitable (generic) conditions on the potential
(like $f_{11,\vzero} \neq 0$) the problem is of comparable difficulty
(cf. \cite{JS,XZ} for Diophantine $\oo$), and reducibility on a
large measure set of values of $\eps$ can be proved. But, in general,
if no condition at all is assumed on the potential (besides analyticity),
further difficulties arise; cf. \cite{GBC,Ch,GGG} for similar situations.
In particular in \cite{GBC} Hill's equation perturbed with a small
quasi-periodic potential was studied under the standard Diophantine
condition, and reducibility for a Cantor set of values of $\eps$ was proved.

More generally one can consider skew-products flows on $\TTT^{d} \times
\SL(n,\RRR)$, for any $n$ (and any $d$). In principle our techniques
apply also in such a case: of course the tree formalism becomes more
involved. Also, less smooth potentials can be considered,
like in \cite{K1,KL,AK}, but in the case of Bryuno vectors analyticity
is likely to be the optimal regularity condition on the potential.

\zerarcounters
\section{Preliminary considerations}
\label{sec:2}

Assume $\la \in [a,b] \subset \RRR \setminus \{0\}$; we shall see later that
the condition $0\notin[a,b]$ can be relaxed (cf. the end  of Section
\ref{sec:6}). Let $A \in \esselle(2,\RRR)$ with imaginary eigenvalues.
Possibly renaming $a$ and $b$ we can assume that the eigenvalues
be $\pm i$. Let $f\!:\TTT^{d}\to\esselle(2,\RRR)$ be real-analytic,
$\oo\in \RRR^{d}$ a real vector, and $\eps$ a real parameter.

Consider the ordinary differential equation
\begin{equation}
\dot x = \left( \la A + \eps f(\oo t) \right) x , 
\label{eq:2.1} \end{equation}
on $\SL(2,\RRR)$.

We can assume that $A$ be of the form
\begin{equation}
A = \left( \begin{matrix} 0 & 1 \\ -1 & 0 \end{matrix} \right) ,
\label{eq:2.2} \end{equation}
and, through a suitable change of coordinates, we obtain
\begin{equation}
D := M A M^{-1} = 
\left( \begin{matrix} i & 0 \\ 0 & -i \end{matrix} \right) ,
\qquad M = \frac{1}{2}
\left( \begin{matrix} 1 & -i \\ 1 & i \end{matrix} \right) ,
\qquad M^{-1} =
\left( \begin{matrix} 1 & 1 \\ i & -i \end{matrix} \right) ,
\label{eq:2.3} \end{equation}
Then, for $z=MxM^{-1}$, we find the equation
\begin{equation}
\dot z = \left( \la D + \eps g(\oo t) \right) z , 
\label{eq:2.4} \end{equation}
with $g=MfM^{-1}$.

Let us introduce some notations.
Given a $2\times2$ matrix $M$, we write
\begin{equation}
M = \left( \begin{matrix} M_{11} & M_{12} \\
M_{21} & M_{22} \end{matrix} \right) ,
\label{eq:2.5} \end{equation}
and we denote by $[A,B]$ the commutator of the two matrices $A$ and $B$.
For $z\in \CCC$ denote by $z^{*}$ the complex conjugate of $z$.
By $\delta_{i,j}$ we denote the Kronecker delta.
We set $\ZZZ_{+}=\{n\in\ZZZ : n \ge 0\}=\NNN \cup \{0\}$, and
for $d \in \NNN$ and $\vzero \in \ZZZ^{d}$, define
$\ZZZ^{d}_{*}=\ZZZ^{d} \setminus\{\vzero\}$.
Given any set $A\subset \RRR$, we denote by $\hbox{meas}(A)$
the Lebesgue measure of $A$.

\begin{lemma}\label{lem:1}
Let $g=MfM^{-1}$, with $f\in\esselle(2,\RRR)$ and $M$
given as in (\ref{eq:2.3}). Then $g\in\esselle(2,\CCC)$,
and one has $g_{11}=g_{22}^{*}$ and $g_{12}=g_{21}^{*}$.
\end{lemma}

\prova The property for $g$ to be traceless follows from the fact
that $\tr (MfM^{-1})= \tr f = 0$. The relations between
the entries of $g$ can be checked by a direct computation:
\begin{eqnarray}
2 g_{11} & = & f_{11} + f_{22} + i \left( f_{12} - f_{21} \right) ,
\nonumber \\
2 g_{12} & = & f_{11} - f_{22} - i \left( f_{12} + f_{21} \right) ,
\nonumber \\
2 g_{21} & = & f_{11} - f_{22} + i \left( f_{12} + f_{21} \right) ,
\nonumber \\
2 g_{22} & = & f_{11} + f_{22} - i \left( f_{12} - f_{21} \right) ,
\label{eq:2.6} \end{eqnarray}
where all entries $f_{ij}$ are real.$\EP$

Define
\begin{eqnarray}
\gotM & := & \left\{ G \in \SL(2,\CCC) : G_{11}=G_{22}^{*},
\quad G_{12}=G_{21}^{*} \right\} , \nonumber \\
\gotm & := & \left\{ g \in \esselle(2,\CCC) : g_{11}=g_{22}^{*},
\quad g_{12}=g_{21}^{*} \right\} .
\label{eq:2.7} \end{eqnarray}
It is easy to see that $\gotM$ is a subgroup, and $\gotm$
is the corresponding Lie algebra.

\begin{lemma}\label{lem:2}
Consider the equation $\dot z = Sz$, with $S=S(t) \in \gotm$
and $z(0) \in \gotM$. Then $z(t) \in \gotM$ for all $t\in\RRR$
for which the solution is defined.
\end{lemma}

\prova Write explicitly the equations for the entries of $z$:
\begin{eqnarray}
\dot z_{11} & = & S_{11} z_{11} + S_{12} z_{21} , \nonumber \\
\dot z_{12} & = & S_{11} z_{12} + S_{12} z_{22} , \nonumber \\
\dot z_{21} & = & S_{21} z_{11} + S_{22} z_{21} =
S_{12}^{*} z_{11} + S_{11}^{*} z_{21} , \nonumber \\
\dot z_{22} & = & S_{21} z_{12} + S_{22} z_{22} =
S_{12}^{*} z_{12} + S_{11}^{*} z_{22} ,
\label{eq:2.8} \end{eqnarray}
so that, by setting $w=(w_{1},w_{2})$, with $w_{1}=z_{11}-z_{22}^{*}$
and $w_{2}=z_{21}-z_{12}^{*}$, one obtains $\dot w = S w$.
If $z(0)\in\gotM$ then $w(0)=0$, so that $w(t)=0$
for all $t\in\RRR$. Moreover, if $\delta(t)=\det z(t)$, one finds
\begin{equation}
\dot \delta = \left( S_{11} + S_{11}^{*} \right)
\left( z_{11} z_{22} - z_{12} z_{21} \right) =
\left( S_{11} + S_{11}^{*} \right) \delta ,
\label{eq:2.9} \end{equation}
where $S_{11} + S_{11}^{*}=S_{11}+S_{22}=\tr S = 0$.
Hence $\delta(t)=\delta(0)=1$. $\EP$

Therefore it is not restrictive to consider the differential equation
\begin{equation}
\dot x = \left( \la A + \eps f(\oo t) \right) x , 
\label{eq:2.10} \end{equation}
on $\gotM$, with
\begin{equation}
A = \left( \begin{matrix} i & 0 \\ 0 & -i \end{matrix} \right) ,
\qquad f \in C^{\om}(\TTT^{d},\gotm) ,
\label{eq:2.11} \end{equation}
and this we shall do henceforth. Write $\la=\la_{0}+\mu$,
and set $x=B(\oo t) y$, with $y$ solution of
\begin{equation}
\dot y = \la_{0} A y , \qquad y(0) = \one ,
\label{eq:2.12} \end{equation}
that is
\begin{equation}
y(t) = 
\left( \begin{matrix} {\rm e}^{i \la_{0} t} & 0 \\ 
0 & {\rm e}^{-i \la_{0} t} \end{matrix} \right) .
\label{eq:2.13} \end{equation}
Then $B = B(\oo t)$ must solve the differential equation
\begin{equation}
\dot  B + \la_{0} [B,A] = \left( \eps f + \mu A \right) B ,
\label{eq:2.14} \end{equation}
and one has $\det B=1$ if $\det x(0) = 1$.

Write
\begin{equation}
B := 1 + \beta , \qquad \beta =
\left( \begin{matrix} a & b \\ c & d \end{matrix} \right) .
\label{eq:2.15} \end{equation}
Then the following result holds.

\begin{lemma}\label{lem:3}
With the notations in (\ref{eq:2.15}) one has $a=d^{*}$, $b=c^{*}$,
and $a + d + ad - bc$ is constant along the flow.
If $\det B =1$ then
\begin{equation}
\tr \beta + \det \beta = a + d + \left( ad - bc \right) = 0 .
\label{eq:2.16} \end{equation}
for all $t\in\RRR$.
\end{lemma}

\prova Since $\gotM$ is a group and $y \in \gotM$,
then $B \in \gotM$ if $x\in\gotM$. More generally,
$\det B(\oo t) = \det B(0)$, which means that
$\det B = 1 + a + d + ad -bc$ is constant along the flow. 
By requiring $\det B=1$ gives (\ref{eq:2.16}).$\EP$

In terms of $\beta$, (\ref{eq:2.14}) becomes

\begin{equation}
\dot \beta + \la_{0} [\beta,A] = \left( \eps f + \mu A \right)
\left( 1 + \beta \right) ,
\label{eq:2.17} \end{equation}
which, written explicitly for the corresponding entries, gives
\begin{eqnarray}
\dot a
& = & \eps f_{11} + i \mu + \eps
\left( f_{11} a + f_{12} c \right) + i \mu \, a , \nonumber \\
\dot b - 2 i \la_{0} b
& = & \eps f_{12} + \eps
\left( f_{11} b + f_{12} d \right) + i \mu \, b , \nonumber \\
\dot c + 2 i \la_{0} c
& = & \eps f_{21} + \eps
\left( f_{21} a + f_{22} c \right) - i \mu \, c , \nonumber \\
\dot d
& = & \eps f_{22} - i \mu + \eps
\left( f_{21} b + f_{22} d \right) - i \mu \, d .
\label{eq:2.18} \end{eqnarray}
If we use that $d=a^{*}$ and $b=c^{*}$, equations (\ref{eq:2.18})
reduce to two independent equations
\begin{eqnarray}
\dot a
& = & \eps f_{11} + i \mu + \eps
\left( f_{11} a + f_{12} c \right) + i \mu \, a , \nonumber \\
\dot c + 2 i \la_{0} c
& = & \eps f_{21} + \eps
\left( f_{21} a + f_{22} c \right) - i \mu \, c ,
\label{eq:2.19} \end{eqnarray}
which is the system the we are going to study.

We can view (\ref{eq:2.19}) as a system of ordinary
differential equations on $\CCC^{2}$.
Such a system admits a first integral, as the following result shows.

\begin{lemma}\label{lem:4}
Given the system (\ref{eq:2.19}), the function
\begin{equation}
H = H(a,c) := a + a^{*} + \left( |a|^{2} - |c|^{2} \right)
\label{eq:2.20} \end{equation}
is a constant of motion, that is $\dot H=0$.
\end{lemma}

\prova Just note that (\ref{eq:2.19}) is a rewriting of
(\ref{eq:2.14}). Lemma \ref{lem:3} shows that $\det B$
is a constant of motion. In terms of $a$ and $c$, this
means that (\ref{eq:2.20}) is conserved along the flow.$\EP$

\zerarcounters
\section{Formal series}
\label{sec:3}

For any function $F$ defined on $\TTT^{d}$ set, formally,
\begin{equation}
F(\pps) = \sum_{\nn \in \ZZZ^{d}}
{\rm e}^{i\nn\cdot\pps} F_{\nn} ,
\label{eq:3.1} \end{equation}
where $\cdot$ denotes the standard inner product in $\RRR^{d}$.
If $F$ is analytic the Fourier coefficients $F_{\nn}$ decay exponentially
at infinity. In particular if $f\in C^{\om}(\TTT^{d},\esselle)$
there exists two constants $F_{0}$ and $\ka_{0}$ such that
$|f_{jj',\nn}| \le F_{0} {\rm e}^{-\ka_{0}|\nn|}$ for $j,j'=1,2$.

Assume that $\oo\in\RRR^{d}$ is a \textit{Bryuno vector}. This means that,
by setting $\al_{n}(\oo)=\inf_{|\nn| \le 2^{n}}|\oo\cdot\nn|$, one has
\begin{equation}
\BB(\oo) := \sum_{n=0}^{\infty} \frac{1}{2^{n}} \log
\frac{1}{\al_{n}(\oo)} < \infty .
\label{eq:3.2} \end{equation}

In terms of the Fourier coefficients $\beta_{\nn}$,
(\ref{eq:2.19}) gives for $\nn\neq \vzero$
\begin{eqnarray}
i\oo\cdot\nn \, a_{\nn} & = & \eps f_{11,\nn} +
\eps \left( f_{11} a + f_{12} c \right)_{\nn} + i \mu \, a_{\nn} ,
\nonumber \\
i \left( \oo\cdot\nn + 2 \la_{0} \right) c_{\nn} & = &
\eps f_{21,\nn} +
\eps \left( f_{21} a + f_{22} c \right)_{\nn} - i \mu \, c_{\nn} ,
\label{eq:3.3} \end{eqnarray}
and for $\nn=\vzero$
\begin{eqnarray}
0 & = & \eps f_{11,\vzero} + i \mu + 
\eps \left( f_{11} a + f_{12} c \right)_{\vzero} + i \mu \, a_{\vzero} ,
\nonumber \\
2 i \la_{0} \, c_{\vzero} & = &
\eps f_{21,\vzero} +
\eps \left( f_{21} a + f_{22} c \right)_{\vzero} - i \mu \, c_{\vzero} .
\label{eq:3.4} \end{eqnarray}
%

\subsection{Recursive equations}

Assume $\la\neq0$. We shall see that $\mu=O(\eps)$, so that
the assumption is satisfied for all $\la\in[a,b]$ if 
$\eps$ is small enough and $0\notin[a,b]$. In fact it would be enough
to require that $\min\{|a|,|b|\}$ be of order $|\eps|^{\sigma}$;
cf. the end of Section \ref{sec:6}.

We can write a formal power series in $\eps$ for $\beta$, by setting
\begin{equation}
\beta = \beta(\oo t) =
\sum_{k=1}^{\infty} \eps^{k} \beta^{(k)}(\oo t) ,
\qquad \beta^{(k)}(\pps) =
\sum_{\nn\in\ZZZ^{d}} {\rm e}^{i\nn\cdot\pps} \beta^{(k)}_{\nn} .
\label{eq:3.5} \end{equation}
The properties $a=d^{*}$ and $b=c^{*}$ imply
$a_{\nn}^{*}=d_{-\nn}$ and $b_{\nn}^{*}=c_{-\nn}$.
In the same way $f\in\gotm$ yields $f_{11,\nn}^{*}=f_{22,-\nn}$,
hence $f_{11,\nn}+f_{11,-\nn}^{*}=0$, and
$f_{12,\nn}^{*}=f_{21,-\nn}$.

If we write also
\begin{equation}
\mu = \sum_{k=1}^{\infty} \eps^{k} \mu^{(k)} ,
\label{eq:3.6} \end{equation}
and we insert (\ref{eq:3.5}) and (\ref{eq:3.6})
into (\ref{eq:3.3}) and (\ref{eq:3.4}) we find
\begin{eqnarray}
a^{(1)}_{\nn} & = & - i \frac{f_{11,\nn}}{\oo\cdot\nn} , \nonumber \\
c^{(1)}_{\nn} & = & - i \frac{f_{21,\nn}}{\oo\cdot\nn + 2\la_{0}} ,
\label{eq:3.7} \end{eqnarray}
for $k=1$ and $\nn\neq \vzero$,
\begin{eqnarray}
\mu^{(1)} & = & i \, f_{11,\vzero} , \nonumber \\
c^{(1)}_{\vzero} & = & - \frac{i}{2 \la_{0}} \, f_{21,\vzero} ,
\label{eq:3.8} \end{eqnarray}
for $k=1$ and $\nn= \vzero$,
\begin{eqnarray}
a^{(k)}_{\nn} & = & - i \frac{1}{\oo\cdot\nn} 
\left( \sum_{\nn_{1}+\nn_{2} = \nn}
\left( f_{11,\nn_{1}} a^{(k-1)}_{\nn_{2}} +
f_{12,\nn_{1}} c^{(k-1)}_{\nn_{2}} \right) 
+ i \sum_{k_{1}+k_{2}=k} \mu^{(k_{1})} a^{(k_{2})}_{\nn} \right) ,
\nonumber \\
c^{(k)}_{\nn} & = & - i \frac{1}{\oo\cdot\nn + 2\la_{0}} 
\left( \sum_{\nn_{1}+\nn_{2} = \nn}
\left( f_{21,\nn_{1}} a^{(k-1)}_{\nn_{2}} +
f_{22,\nn_{1}} c^{(k-1)}_{\nn_{2}} \right) 
- i \sum_{k_{1}+k_{2}=k} \mu^{(k_{1})} c^{(k_{2})}_{\nn} \right) ,
\label{eq:3.9} \end{eqnarray}
for $k\ge2$ and $\nn\neq \vzero$, and
\begin{eqnarray}
\mu^{(k)} & = & i \left( \sum_{\nn_{1}+\nn_{2} = \vzero}
\left( f_{11,\nn_{1}} a^{(k-1)}_{\nn_{2}} +
f_{12,\nn_{1}} c^{(k-1)}_{\nn_{2}} \right) +
i \sum_{k_{1}+k_{2}=k} \mu^{(k_{1})} a^{(k_{2})}_{\vzero} \right) ,
\nonumber \\
c^{(k)}_{\vzero} & = & - \frac{i}{2 \la_{0}}
\left( \sum_{\nn_{1}+\nn_{2} = \vzero}
\left( f_{21,\nn_{1}} a^{(k-1)}_{\nn_{2}} +
f_{22,\nn_{1}} c^{(k-1)}_{\nn_{2}} \right) -
i \sum_{k_{1}+k_{2}=k} \mu^{(k_{1})} c^{(k_{2})}_{\vzero} \right) ,
\label{eq:3.10} \end{eqnarray}
for $k\ge2$ and $\nn=\vzero$.

\begin{lemma}\label{lem:5}
Let $\mu^{(k)}$ and $c^{(k)}_{\vzero}$ be fixed for all $k\ge 1$
according to (\ref{eq:3.8}) and (\ref{eq:3.10}). Then
there are formal power series (\ref{eq:3.5}) and (\ref{eq:3.6})
for $\beta$ and $\mu$, respectively,
recursively determined from (\ref{eq:3.7}) to (\ref{eq:3.10}),
which solve order by order equations (\ref{eq:2.19}).
The constants $a^{(k)}_{\vzero}$ can be arbitrarily fixed.
\end{lemma}

We omit the easy proof, which can be obtained also as a byproduct
of the forthcoming analysis. Therefore, the formal solubility of
the equations (\ref{eq:2.19}) requires that $\mu^{(k)}$
and $c^{(k)}_{\vzero}$ be fixed to all orders $k \ge1$,
while all coefficients $a^{(k)}_{\vzero}$ are left undetermined.
We can fix the latter by requiring (\ref{eq:2.16}).

\begin{lemma}\label{lem:6}
If we fix
\begin{equation}
a^{(1)}_{\vzero} = 0
\label{eq:3.11} \end{equation}
for $k=1$, and
\begin{equation}
a^{(k)}_{\vzero} = - \frac{1}{2} \sum_{k_{1}+k_{2}=k}
\sum_{\nn\in\RRR^{d}}
\left( a^{(k_{1})}_{\nn} a^{(k_{2})*}_{\nn} -
c^{(k_{1})}_{\nn} c^{(k_{2})*}_{\nn} \right) ,
\label{eq:3.12} \end{equation}
for $k\ge 2$, then
\begin{equation}
H^{(k)} := a^{(k)} + a^{(k)*} + \sum_{k_{1}+k_{2}=k} \left(
a^{(k_{1})} a^{(k_{2})*} - c^{(k_{1})} c^{(k_{2})*} \right) = 0,
\label{eq:3.13} \end{equation}
for all $k \in \NNN$.
\end{lemma}

\prova By Lemma \ref{lem:4} to all orders $k\ge 1$ the function $H^{(k)}$
is formally a constant, so that $H^{(k)}_{\nn}=0$ for all $k\ge 1$
and all $\nn\neq\vzero$, while $H^{(1)}_{\vzero}=a^{(1)}_{\vzero}+
a^{(1)*}_{\vzero}$ and
\begin{equation}
H^{(k)}_{\vzero} = a^{(k)}_{\vzero} + a^{(k)*}_{\vzero} +
\sum_{k_{1}+k_{2}=k} \left( a^{(k_{1})} a^{(k_{2})*} -
c^{(k_{1})} c^{(k_{2})*} \right)_{\vzero}
\label{eq:3.14} \end{equation}
for $k\ge 2$ are constants. If we fix $a^{(k)}_{\vzero}$ recursively
according to (\ref{eq:3.11}) and (\ref{eq:3.12}),
then $H^{(k)}_{\vzero}=0$, so that (\ref{eq:3.13}) follows.$\EP$

The recursive equations (\ref{eq:3.7}) to (\ref{eq:3.12})
can be graphically represented in terms of linear trees as follows.

\begin{figure}[ht]
\vskip.3truecm
\centering
\ins{75pt}{-2pt}{(a)}
\ins{130pt}{13pt}{$(k)$}
\ins{110pt}{-8pt}{$j,\vzero$}
\ins{180pt}{-2pt}{(b)}
\ins{238pt}{13pt}{$(1)$}
\ins{218pt}{-8pt}{$j,\nn$}
\ins{290pt}{-2pt}{(c)}
\ins{347.pt}{13pt}{$(k)$}
\ins{328pt}{-8pt}{$j,\nn$}
\includegraphics{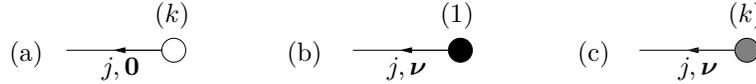}
\caption{Graph elements representing (a) $u^{(k)}_{j,\vzero}$
for $j=1,2$ and $\mu^{(k)}$ for $j=3$, (b) $u^{(1)}_{j,\nn}$,
$\nn\neq\vzero$, and (c) $u^{(k)}_{j,\nn}$. Only in (a) one can
have $j=3$, otherwise $j=1,2$. For $\nn=0$ the latter graph
reduces to the first graph, while for $k=1$ and $\nn\ne\vzero$
it reduces to the second graph.}
\label{fig:3.1}
\end{figure}

Call $u=(u_{1},u_{2})=(a,c)$. Set also $u^{(k)}_{1,\nn}=a^{(k)}_{\nn}$
and $u^{(k)}_{2,\nn}=c^{(k)}_{\nn}$, and represent $u_{j,\nn}^{(k)}$
as a line carrying the labels $j\in\{1,2\}$ and $\nn\in\ZZZ^{d}$
exiting from a bullet carrying the label $k$, with $k\in\NNN$.
We call $k,j,\nn$ the \textit{order} label, the \textit{component} label
and the \textit{momentum} label, respectively.
We colour the bullet with white if $\nn=\vzero$ and
with grey if $\nn\neq\vzero$; in the latter case, for $k=1$
we draw the bullet as a black bullet instead of
a grey one; cf. Figure \ref{fig:3.1}.
We call \textit{graph elements} the graphs which are drawn this way.
We represent also $\mu^{(k)}$ by a graph element,
by using the same graph for $u^{(k)}_{\vzero}$ except that $j=3$,
i.e. we set $\mu^{(k)}=u^{(k)}_{3,\vzero}$.

\begin{figure}[ht]
\vskip.3truecm
\centering
\ins{72pt}{-8pt}{$j,\nn$}
\ins{92pt}{13pt}{$(k)$}
\ins{128pt}{-5pt}{$=$}
\ins{231pt}{13pt}{$(k\!\!-\!\!1)$}
\ins{168pt}{-8pt}{$j,\nn$}
\ins{196pt}{-11pt}{$\nn_{1}$}
\ins{211pt}{-8pt}{$j_{2},\nn_{2}$}
\ins{277pt}{-2pt}{$+$}
\ins{338pt}{13pt}{$(k_{1})$}
\ins{317pt}{-8pt}{$j,\nn$}
\ins{384pt}{13pt}{$(k_{2})$}
\ins{363pt}{-8pt}{$j,\nn$}
\includegraphics{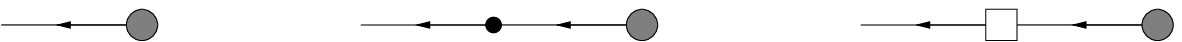}
\caption{Graphical representation of (\ref{eq:3.9}),
expressing the coefficient $u_{j,\nn}^{(k)}$ for $k \ge 2$, $j=1,2$,
and $\nn\neq\vzero$ in terms of the coefficients $u_{j',\nn'}^{(k')}$,
with $k'<k$. In the first graph one has the constraint $\nn=
\nn_{1}+\nn_{2}$, while in second graph one has
the constraint $k=k_{1}+k_{2}$.}
\label{fig:3.2}
\end{figure}

Then equations (\ref{eq:3.9}) can be represented as shown in
Figure \ref{fig:3.2}, provided we give some rules
in order to associate with the graphs suitable numerical values.

In the two graphs on the right hand side of Figure \ref{fig:3.2}
there are two lines $\ell_{1}$ and $\ell_{2}$,
with labels $(j_{\ell_{1}},\nn_{\ell_{1}})=(j,\nn)$
and $(j_{\ell_{2}},\nn_{\ell_{2}})=(j_{2},\nn_{2})$, respectively.
In the first graph we associate a \textit{mode} label
$\nn_{v}=\nn_{1}\in\ZZZ^{d}$ and a \textit{node factor}
$F_{v}=f_{jj_{2},\nn_{1}}$ with the black point $v$ between the two lines.
In the second graph we associate a \textit{mode} label
$\nn_{v}=\nn_{1}= \vzero$, an \textit{order} label $k_{v}=k_{1}$
and a \textit{node factor} $F_{v}=(-1)^{j+1} i \mu^{(k_{1})}
\delta_{j,j_{2}}$ with the white square $v$ between the two lines.
In both graphs we have the constraint $\nn=\nn_{1}+\nn_{2}$,
which fixes $\nn_{2}=\nn$ in the second graph.
With the line $\ell_{1}$ we associate a \textit{propagator}
$g_{\ell_{1}}$, such that $g_{\ell_{1}}=1/i\oo\cdot\nn$
if $j=1$ and $g_{\ell_{1}}=1/i(\oo\cdot\nn+2\la_{0})$
if $j=2$, -- note that in both graphs one has $\nn\neq\vzero$.
Finally the line $\ell_{2}$ together with the grey bullet which it
comes out from forms a graph element as shown in Figure \ref{fig:3.1}(c),
so that it represents $u^{(k_{2})}_{j_{2},\nn_{2}}$, with $k_{2}=k-1$
in the first graph and $(j_{2},\nn_{2})=(j,\nn)$ in the second one.

To obtain $u^{(k)}_{j,\nn}$, with $\nn\neq\vzero$, one has to sum over
all labels the products of the propagator $g_{\ell_{1}}$ times the
node factor $F_{v}$ times the coefficient $u^{(k_{2})}_{j_{2},\nn_{2}}$
represented by the graph element attached either to the black point
or to the white square, with the constraint that the labels $j,\nn,k$
are kept fixed. The quantity that one obtains this way is just the right
hand side of equations (\ref{eq:3.9}). Of course $j=1$ means that the
corresponding graphs represent contributions to $a^{(k)}_{\nn}$,
and $j=2$ means that they represent contributions to $c^{(k)}_{\nn}$.

\begin{figure}[ht]
\vskip.3truecm
\centering
\ins{70pt}{-8pt}{$j,\vzero$}
\ins{92pt}{13pt}{$(k)$}
\ins{130pt}{-3pt}{$=$}
\ins{231pt}{13pt}{$(k\!\!-\!\!1)$}
\ins{172pt}{-8pt}{$j,\vzero$}
\ins{197pt}{-11pt}{$\nn_{1}$}
\ins{213pt}{-8pt}{$j_{2},\nn_{2}$}
\ins{280pt}{-2pt}{$+$}
\ins{339pt}{13pt}{$(k_{1})$}
\ins{318pt}{-8pt}{$j,\vzero$}
\ins{384pt}{13pt}{$(k_{2})$}
\ins{361pt}{-8pt}{$j,\vzero$}
\includegraphics{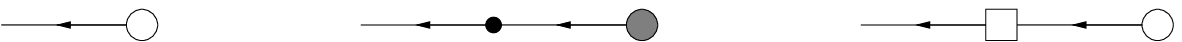}
\caption{Graphical representation of (\ref{eq:3.10}),
expressing the constants $\mu^{(k)}$ (if $j=3$)
and $c^{(k)}_{\vzero}$ (if $j=2$) for $k \ge 2$
in terms of the coefficients $u_{j',\nn'}^{(k')}$, with $k'<k$.
In the first graph one has the constraint $\vzero=
\nn_{1}+\nn_{2}$, while in second graph one has
the constraint $k=k_{1}+k_{2}$.}
\label{fig:3.3}
\end{figure}

Analogously we can represent graphically (\ref{eq:3.10}) as in
Figure \ref{fig:3.3}. The difference with respect to
Figure \ref{fig:3.2} is that now $\nn=\vzero$, and $j\in\{2,3\}$.
For $j=3$ we obtain a contribution to $\mu^{(k)}$, whereas for $j=2$
we have a contribution to $c^{(k)}_{\vzero}$. The quantities to be
associated with the black points, the white bullets, the white squares
and the graph elements are the same as defined in the case
of Figure \ref{fig:3.2}. With the line $\ell_{1}$
we associate a propagator $g_{\ell_{1}}$, such that
$g_{\ell_{1}}=i$ if $j=3$ and $g_{\ell_{1}}=-i/2\la_{0}$ if $j=2$.

\begin{figure}[ht]
\vskip.3truecm
\centering
\ins{145pt}{-29pt}{$1,\vzero$}
\ins{168pt}{-08pt}{$(k)$}
\ins{202pt}{-24pt}{$=$}
\ins{250pt}{-29pt}{$1,\vzero$}
\ins{273pt}{-29pt}{$\vzero$}
\ins{299pt}{-15pt}{$j_{1},\nn_{1}$}
\ins{309pt}{12pt}{$(k_{1})$}
\ins{280pt}{-39pt}{$j_{2},\nn_{2}$}
\ins{309pt}{-26pt}{$(k_{2})$}
\includegraphics{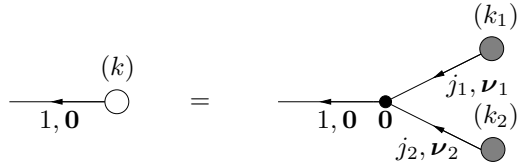}
\caption{Graphical representation of (\ref{eq:3.12}),
expressing the constant $a^{(k)}_{\vzero}$ (hence $j=1$) for $k \ge 2$
in terms of the coefficients $u_{j',\nn'}^{(k')}$, with $k'<k$.
One has the constraints $\vzero=\nn_{1}+\nn_{2}$,
$k=k_{1}+k_{2}$, and $j_{1}=j_{2} \in \{1,2\}$.}
\label{fig:3.4}
\end{figure}

Finally, also $a_{\vzero}^{(k)}$ can be graphically represented
from equation (\ref{eq:3.12}) in terms of the coefficients with
lower order; cf. Figure \ref{fig:3.4}. In such a case, in the graph
on the right hand side, the line $\ell_{1}$ which carries the labels
$(j_{\ell_{1}},\nn_{\ell_{1}})=(1,\vzero)$ has propagator
$g_{\ell_{1}}=1/2$, and comes out from a white bullet $v$
with two entering lines carrying labels $(j_{\ell_{2}},\nn_{\ell_{2}})=
(j_{1},\nn_{1})$ and $(j_{\ell_{3}},\nn_{\ell_{3}})=(j_{2},\nn_{2})$,
with the constraints $j_{1}=j_{2}$ and $\nn_{1}+\nn_{2}=\vzero$.
The node factor is $F_{v}=(-1)^{j}$.

\subsection{Linear trees}

We can iterate the graphical construction given in
Figures \ref{fig:3.2}, \ref{fig:3.3} and \ref{fig:3.4} by developing
further the graph elements on the right hand side according to
same figures. At the end we obtain that
$u^{(k)}_{j,\nn}$, $\nn\neq0$, $\mu^{(k)}$ and $c^{(k)}_{\vzero}$
can all be expressed in terms of linear trees (or chains),
which are constructed as follows.

A \textit{tree} is a collection of \textit{points} and \textit{lines}
connecting them, such that all lines are oriented toward a unique point,
with the property that only one line enters such a point.
The latter is called the \textit{root} of the tree,
and the line entering the root is called the \textit{root line}.
By construction any point
different from the root has one and only one line coming out from it,
called the \textit{exiting line} of the point. A \textit{linear tree}
is a tree such that each point has only one line
going into it, called the \textit{entering line} of the point,
except one which has no entering line at all. The latter is called
the \textit{endpoint} of the tree. All the points except the root and
the endpoint are called the \textit{nodes} of the tree.

Denote by  $V(\theta)$ and $L(\theta)$ the set of nodes and
the set of lines, respectively, in the tree $\theta$.
One has $|L(\theta)|=|V(\theta)|+1$.
Sometimes it can be convenient to denote by $P(\theta)$ the set
of nodes plus the endpoint of $\theta$.

We can number the lines and nodes as $\ell_{1},\ldots,\ell_{N}$,
and $v_{1},\ldots,v_{N-1}$, with $N=|L(\theta)| \ge 1$, in such a way
that $\ell_{N}$ connects the endpoint $v_{N}$ to
the node $v_{N-1}$ (the \textit{first node}),
each line $\ell_{k}$, $k=2,\ldots,N-1$, connects the node $v_{k}$
to the node $v_{k-1}$, and $\ell_{1}$ connects the node $v_{1}$
(the \textit{last node}) to the root.

\begin{figure}[ht]
\centering
\vskip.3truecm
\ins{70pt}{-8pt}{$j,\nn$}
\ins{94pt}{-10pt}{$\nn_{1}$}
\ins{110pt}{-8pt}{$j_{2},\nn_{2}$}
\ins{132pt}{13pt}{$(k_{1})$}
\ins{151pt}{-8pt}{$j_{2},\nn_{2}$}
\ins{180pt}{-10pt}{$\nn_{3}$}
\ins{194pt}{-8pt}{$j_{4},\nn_{4}$}
\ins{221pt}{-10pt}{$\nn_{5}$}
\ins{237pt}{-8pt}{$j_{6},\nn_{6}$}
\ins{259pt}{13pt}{$(k_{2})$}
\ins{278pt}{-8pt}{$j_{6},\nn_{6}$}
\ins{305pt}{-10pt}{$\nn_{7}$}
\ins{320pt}{-8pt}{$j_{8},\nn_{8}$}
\ins{347pt}{-10pt}{$\nn_{9}$}
\ins{360pt}{-8pt}{$j_{1\!0},\nn_{1\!0}$}
\includegraphics{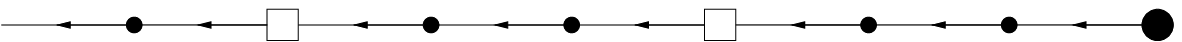}
\caption{An example of tree of order $k$ with $7$ nodes and $8$ lines,
and with an endpoint which is a black bullet.
One has the constraints $k=6+k_{1}+k_{2}$,
$\nn=\nn_{1}+\nn_{2}$, $\nn_{2}=\nn_{3}+\nn_{4}$,
$\nn_{4}=\nn_{5}+\nn_{6}$, $\nn_{6}=\nn_{7}+\nn_{8}$,
$\nn_{8}=\nn_{9}+\nn_{10}$. The constraint that
the lines connected to the white squares carry the same
component and momentum labels has been taken into account explicitly.
The order labels of the black points and of the black bullet
are not shown, as they are necessarily $1$. Also the mode label
of the black bullet is not shown, as it is necessarily $\nn_{1\!0}$}
\label{fig:3.5}
\end{figure}

\begin{figure}[ht]
\centering
\vskip.3truecm
\ins{132pt}{-8pt}{$j,\nn$}
\ins{155pt}{-10pt}{$\nn_{1}$}
\ins{170pt}{-8pt}{$j_{2},\nn_{2}$}
\ins{196pt}{-10pt}{$\nn_{3}$}
\ins{211pt}{-8pt}{$j_{4},\nn_{4}$}
\ins{234pt}{13pt}{$(k_{1})$}
\ins{252pt}{-8pt}{$j_{4},\nn_{4}$}
\ins{281pt}{-10pt}{$\nn_{5}$}
\ins{298pt}{-8pt}{$j_{6},\vzero$}
\ins{320pt}{13pt}{$(k_{2})$}
\includegraphics{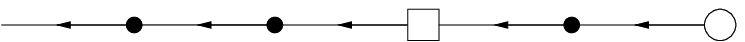}
\caption{An example of tree of order $k$ with $4$ nodes and $5$ lines,
and with an endpoint which is a white bullet.
One has the constraints $k=3+k_{1}+k_{2}$,
$\nn=\nn_{1}+\nn_{2}$, $\nn_{2}=\nn_{3}+\nn_{4}$,
$\nn_{4}=\nn_{5}+\nn_{6}$, with $\nn_{6}=\vzero$.
The constraint that the lines connected to the white square
carry the same component and momentum
labels has been taken into account explicitly.
The order labels of the black points
are not shown, as they are necessarily $1$.}
\label{fig:3.6}
\end{figure}

A node $v$ can be either a \textit{black point} or
a \textit{white square}: in the latter case one must have $\nn_{v}=\vzero$.
The endpoint of the tree can be either a \textit{white bullet}
or a \textit{black bullet}: the line $\ell$ coming out from the endpoint
carries a momentum $\nn_{\ell}=\vzero$ in the first case and a
momentum $\nn_{\ell}\neq\vzero$ in the second one. Examples of trees
are depicted in Figure \ref{fig:3.5} and \ref{fig:3.6}.

With each node $v$ which is a black point we associate
an \textit{order} label $k_{v}=1$ and
a \textit{mode} label $\nn_{v}\in\ZZZ^{d}$, and with each node which is
a white square we associate an \textit{order} label $k_{v}\in\NNN$
and a \textit{mode} label $\nn_{v}=\vzero$.
If the endpoint $v$ is a black bullet
we associate with it an \textit{order} label $k_{v}=1$ and
a \textit{mode} label $\nn_{v}\in\ZZZ^{d}$,
if it is a white bullet we associate with it
an \textit{order} label $k_{v} \in \NNN$ and
a \textit{mode} label $\nn_{v}=\vzero$.
With each line $\ell$ we associate a \textit{component} label
$j_{\ell}\in\{1,2,3\}$ and a \textit{momentum} $\nn_{\ell}\in\ZZZ^{d}$. 
For each node we have the \textit{conservation law}
that the momentum of the exiting line equals the sum of the
mode of the node plus the momentum of the entering line

As the tree is linear, for each node $v$ there are only
one line $\ell=\ell_{v}$ which comes out from it and only
one line $\ell_{v}'$ which enters it. If $v=v_{k}$ this means
that $\ell_{v}=\ell_{k-1}$ and $\ell_{v}'=\ell_{k}$.
With these notations, the conservation law reads
$\nn_{\ell_{v}}=\nn_{v}+\nn_{\ell_{v}'}$.

Once all labels have been assigned,
we associate with each node $v$ the \textit{node factor}
\begin{equation}
F_{v} := \begin{cases}
f_{j_{\ell_{v}},j_{\ell_{v}'},\nn_{v}}
\left( 1 - \delta_{j_{\ell_{v}},3} \right) +
f_{1,j_{\ell_{v}'},\nn_{v}} \delta_{j_{\ell_{v}},3} ,
& \hbox{ $v$ is a black point} , \\
(-1)^{j_{\ell_{v}}+1} i \mu^{(k_{v})} \delta_{\nn_{v},\vzero}
\delta_{j_{\ell_{v}},j_{\ell_{v}'}} ,
& \hbox{ $v$ is a white square} , \\
\end{cases}
\label{eq:3.15} \end{equation}
with the endpoint $v$ the \textit{endpoint factor}
\begin{equation}
F_{v} := \begin{cases}
f_{j_{\ell_{v}}1,\nn_{\ell_{v}}} ,
& \hbox{ $v$ is a black bullet} , \\
u_{j_{\ell_{v}},\vzero}^{(k_{v})} ,
& \hbox{ $v$ is a white bullet} ,
\end{cases}
\label{eq:3.16} \end{equation}
and with each line $\ell$ coming out from a node the \textit{propagator}
\begin{equation}
g_{\ell} := \begin{cases}
-i/\oo\cdot\nn_{\ell} ,
& \nn_{\ell}\neq \vzero , \; j_{\ell}=1 , \\
-i/(\oo\cdot\nn_{\ell}+2\la_{0}) ,
& \nn_{\ell}\neq \vzero , \; j_{\ell}=2 , \\
-i/2\la_{0} , & \nn_{\ell}= \vzero , \; j_{\ell}=2 , \\
i , & \nn_{\ell}= \vzero , \; j_{\ell}=3 ,
\end{cases}
\label{eq:3.17} \end{equation}
and with the line $\ell$ coming out from the
endpoint the \textit{propagator}
\begin{equation}
g_{\ell} := \begin{cases}
-i/\oo\cdot\nn_{\ell} ,
& \nn_{\ell}\neq \vzero , \; j_{\ell}=1 , \\
-i/(\oo\cdot\nn_{\ell}+2\la_{0}) ,
& \nn_{\ell}\neq \vzero , \; j_{\ell}=2 , \\
1 , & \nn_{\ell}= \vzero , \; j_{\ell}=1,2 ,
\end{cases}
\label{eq:3.18} \end{equation}
The propagators (\ref{eq:3.17}) and (\ref{eq:3.18})
are equal as far as $\nn_{\ell}\neq\vzero$, but they are different
when $\nn_{\ell}=\vzero$. 

One has the further constraints
that one can have $\nn_{\ell}=\vzero$ in (\ref{eq:3.17}) only if
$\ell$ is the root line, and $j_{\ell_{v}}=3$ in (\ref{eq:3.15})
again only if $\ell_{v}$ is the root line.
In particular the only lines which can have vanishing momentum
are the root line and the line coming out from the endpoint,
and the only line which can have component label $j=3$
is the root line. Finally if $|P(\theta)|=1$ then the endpoint
of $\theta$ has to be a black bullet.
Define $\Theta^{0}_{k,j,\nn}$ the set of linear trees with
labels $j,\nn$ associated with the root line, and
with $\sum_{v\in V(\theta)} k_{v} = k$.

\begin{lemma}\label{lem:7}
Let $\oo$ be a Bryuno vector. One has
\begin{eqnarray}
& & u^{(k)}_{j,\nn} = \sum_{\theta\in\Theta^{0}_{k,j,\nn}}
\Val (\theta) , \qquad k \ge 1, \; \nn \neq \vzero , \;
j=1,2, \nonumber \\
& & \mu^{(k)} = \sum_{\theta\in\Theta^{0}_{k,3,\vzero}} \Val (\theta) ,
\qquad c^{(k)}_{\vzero} =
\sum_{\theta\in\Theta^{0}_{k,2,\vzero}} \Val (\theta) , \qquad k \ge 1 ,
\label{eq:3.19} \end{eqnarray}
where the \textit{tree value} $\Val(\theta)$ is given by
\begin{equation}
\Val(\theta) = \Big( \prod_{\ell\in L(\theta)} g_{\ell} \Big)
\Big( \prod_{v \in P(\theta)} F_{v} \Big) ,
\label{eq:3.20} \end{equation}
with the propagators $g_{\ell}$ defined by (\ref{eq:3.17}) and
(\ref{eq:3.18}), and the factors $F_{v}$ defined by
(\ref{eq:3.15}) and (\ref{eq:3.16}). One has
$\mu^{(k)} \in \RRR$ for all $k \ge 1$.
\end{lemma}

\prova The only non-trivial statement is that
$\mu^{(k)}$ is real, -- the other assertions
can be easily derived from the discussion above
(or can be proved by induction on $k$).

We prove that $\mu^{(k)} \in \RRR$ by induction.
One has $\mu^{(1)} \in \RRR$ because $\mu^{(1)}=i f_{11,\vzero}$,
and $f_{11,\vzero}$ is purely imaginary.

If $k\ge 2$, for each tree $\theta\in \Theta^{0}_{k,3,\vzero}$
we distinguish three cases: (a) the endnode of $\theta$
is a black bullet, (b) the endnode is a white bullet
and the line coming out from it carries a label $j=1$, and
(c) the endnode is a white bullet and the line coming out
from it carries a label $j=2$.

We discuss first case (a). Given $\theta$ we consider the
tree $\tau=\tau(\theta)$ obtained as follows.
First, detach the root line from the
last node and attach it to the endnode, and change
the orientation of all lines; then the last node of $\theta$
becomes the endnode of $\tau$ (graphically it is
transformed from a black point into a black bullet) and vice versa.
Second, change the sign of all the mode labels.

Of course we can write $ \mu^{(k)} =
\sum_{\theta \in \Theta^{0}_{k,3,\vzero}} \Val(\tau(\theta))$.
If we compare $\tau(\theta)$ with $\theta$ we see that
the propagators are not changed, because the
sum of all the mode labels is zero, i.e. $\sum_{v\in P(\theta)}
\nn_{v}=\vzero$. The node factors corresponding to white squares
$v$ are not changed (they remain $\pm i \mu^{(k_{v})}$), while
the node factors corresponding to black points are changed from
$f_{j_{\ell_{v}}j_{\ell_{v}}',\nn_{v}}$ into
$f_{j_{\ell_{v}}'j_{\ell_{v}},-\nn_{v}}$. The same happens to
the endnode factor, which becomes $f_{1j_{\ell_{v}},-\nn_{v}}$.
Recall that one has $f_{12,-\nn}^{*}=f_{21,\nn}$, $f_{21,-\nn}^{*}=
f_{12,\nn}$, and $f_{11,-\nn}^{*}=-f_{11,\nn}$; moreover $g_{\ell}^{*}=-
g_{\ell}$, as it follows from (\ref{eq:3.17}) and (\ref{eq:3.18}),
and $\mu^{(k_{v})*}=\mu^{(k_{v})}$ by the inductive assumption.

Then if we compute $\Val^{*}(\tau)$ we obtain
$\Val^{*}(\tau)=(-1)^{|L(\theta)|}(-1)^{|J(\theta)|}\Val(\theta)$,
where $L(\theta)$ is the set of lines in $\theta$, and $J(\theta)$
is the set of $v\in P(\theta)$ with $j_{\ell_{v}}=j_{\ell_{v}'}$
(we set $j_{\ell_{v}}'=1$ if $v$ is the endnode), hence
including the white squares. It is immediate to realize that
$|P(\theta) \setminus J(\theta)|$ is even, so that
$|J(\theta)|$ has the same parity as $|P(\theta)|$.
As $|P(\theta)|=|L(\theta)|$ this yields
$\Val^{*}(\tau)= \Val(\theta)$.

In case (b) we can write $\Val(\theta)=\Val(\theta_{1})
\, a^{(k_{1})}_{\vzero}$ for suitable $\theta_{1}$ and $k_{1}$,
with $a^{(k_{1})}_{\vzero}$ real by (\ref{eq:3.12}). More precisely
$\theta_{1}$ is the tree of order $k-k_{1}$ obtained from $\theta$
by detaching the graph element representing $a^{(k_{1})}_{\vzero}$
and replacing the first node with an endpoint. Then we can construct
a tree $\tau_{1}=\tau(\theta_{1})$, and reason for $\theta_{1}$
as done for $\theta$ in case (a).
The same conclusions hold, in particular one finds
$\Val^{*}(\tau_{1})=\Val(\theta_{1})$.

Finally in case (c) we can write $\Val(\theta)=\Val(\theta_{1})
\, c^{(k_{1})}_{\vzero}$ for suitable $\theta_{1}$ and $k_{1}$,
and develop $c^{(k_{1})}_{\vzero}$ in terms of trees
(according with a procedure which will be extensively used
in the following), and so on, until we reach a tree which
belongs to case (a) or case (b), up to the fact it can contain
lines $\ell$ with $\nn_{\ell}=\vzero$ and $g_{\ell}=-i/2\la_{0}$;
see (\ref{eq:3.17}) and (\ref{eq:3.18}). Therefore we can reason as
in the previous cases (a) and (b).

By putting together all the cases, at the end we
obtain $\mu^{(k)}=\mu^{(k)*}$.$\EP$

Note that the set $\Theta^{0}_{k,1,\vzero}$ does not appear
in (\ref{eq:3.19}). This is necessary as the map $\theta\to
\Val(\theta)$ is not defined for $\theta\in\Theta^{0}_{k,1,\vzero}$;
see (\ref{eq:3.17}). In fact, $a^{(k)}_{\vzero}$ cannot be represented
as a sum of values of linear trees, but still we can write
for $k \ge 2$ (and setting $a^{(1)}_{\vzero}=0$)
\begin{equation}
a^{(k)}_{\vzero} = \frac{1}{2} \sum_{k_{1}+k_{2}=k}
\sum_{\nn\in\ZZZ^{d}} \sum_{j=1,2} (-1)^{j}
{\mathop{\sum}_{\theta_{1}\in
\Theta^{0}_{k_{1},1,\nn}}}^{\hskip-.5truecm'}
\Val(\theta_{1})
{\mathop{\sum}_{\theta_{2}\in
\Theta^{0}_{k_{2},1,\nn}}}^{\hskip-.5truecm'}
\Val^{*}(\theta_{2}) ,
\label{eq:3.21} \end{equation}
where $'$ means that we must interpret
\begin{equation}
{\mathop{\sum}_{\theta\in\Theta^{0}_{k,1,\vzero}}}^{\hskip-.3truecm'}
\Val(\theta) := a^{(k)}_{\vzero} .
\label{eq:3.22} \end{equation}
Hence also $a^{(k)}_{\vzero}$
can be expressed in terms of linear trees.

\subsection{Nonlinear trees}

Each node represented by a white square can be further expanded in terms
of trees as follows. First replace the white square $v$ with a
black point and attach to the latter a further graph element representing
$\mu^{(k_{v})}$, if $k_{v}$ is the order label of $v$
(cf. Figure \ref{fig:3.5}), hence the graph element is expressed
in terms of trees according to the first graph in Figure \ref{fig:3.3}.
With the new node $v$, represented by a black point, we associate
a mode label $\nn_{v}=\vzero$ and an order label $k_{v}=0$.

\begin{figure}[ht]
\centering
\ins{129pt}{-36pt}{$j,\nn$}
\ins{149pt}{-14pt}{$(k_{v})$}
\ins{172pt}{-36pt}{$j,\nn$}
\ins{226pt}{-33pt}{$\Longrightarrow$}
\ins{272pt}{-36pt}{$j,\nn$}
\ins{296.5pt}{-38pt}{$\vzero$}
\ins{312pt}{-36pt}{$j,\nn$}
\ins{326pt}{14pt}{$(k_{v})$}
\ins{319pt}{-17pt}{$3,\vzero$}
\includegraphics{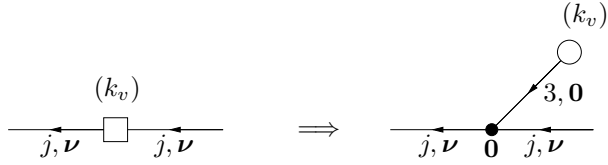}
\vskip.5truecm
\caption{The quantity $\mu^{(k_{v})}$ appearing in the node factor
associated with $v$ can be expressed according to (\ref{eq:3.10}).
This can be interpreted graphically by replacing the white square
as shown in the figure: the graph element entering the node $v$
represent $\mu^{(k_{v})}$, and it can be further developed in terms
of trees according to Figure \ref{fig:3.3}.}
\label{fig:3.7}
\end{figure}

In the same way also the endpoints which are drawn as white bullets
can be expanded according to the second graph in Figure \ref{fig:3.3}
if the exiting line carries a component label $j=2$ and
according to the graph in Figure \ref{fig:3.4}
if the exiting line carries a component label $j=1$.

Of course if we do this, then nonlinear trees appear.
Nonlinear trees are partially ordered sets of points
and lines connecting them, and not totally ordered sets,
such as linear trees are. The advantage
of this procedure, however, is that at the end, the
trees have only endpoints with order $1$ and all the node factors
are quantities fixed (and not to be determined iteratively).
The new trees can have also nodes with two entering lines.
If we denote by $p_{v}$ the \textit{branching number} of the point $v$,
that is the number of lines entering $v$, then $p_{v}=1,2$
if $v$ is a node, while $p_{v}=0$ is $v$ is an endpoint.

\begin{figure}[ht]
\vskip.3truecm
\centering
\ins{28pt}{-13pt}{(i)}
\ins{58pt}{-22pt}{$1,\vzero$}
\ins{100pt}{-11pt}{$j,\nn$}
\ins{85pt}{-32pt}{$j,\nn,*$}
\ins{82.5pt}{-22.5pt}{$\vzero$}
\ins{170pt}{-13pt}{(ii)}
\ins{200pt}{-21pt}{$j,\nn$}
\ins{245pt}{-12pt}{$3,\vzero$}
\ins{238pt}{-32pt}{$j,\nn$}
\ins{226pt}{-22.5pt}{$\vzero$}
\ins{308pt}{-13pt}{(iii)}
\ins{345pt}{-21pt}{$j,\nn$}
\ins{387pt}{-12pt}{$j,\nn$}
\ins{380pt}{-33pt}{$3,\vzero$}
\ins{369.5pt}{-22.5pt}{$\vzero$}
\includegraphics{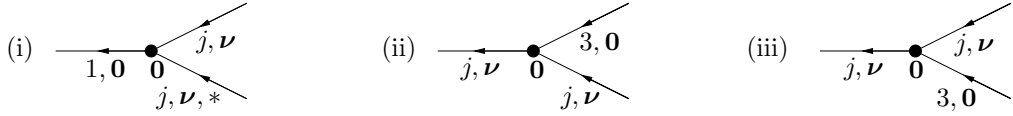}
\caption{Nodes with branching number $2$. The corresponding
node factors are defined in (\ref{eq:3.23}). The entering lines
are assumed to come out from other nodes or from endpoints,
and the exiting line either enters another node or is the root line.
The label $*$ on the lower entering line of the first graph
means that one has to take the complex conjugate of the value
of the subtree with that root line.}
\label{fig:3.8}
\end{figure}

A node $v$ with $p_{v}=2$ has the following properties.
Denote by $\ell_{0}$ the exiting line of $v$, 
and by $\ell_{1}$ and $\ell_{2}$ the entering lines of $v$.
Then either (i) $j_{\ell_{0}}=1$, $\nn_{\ell_{0}}=\vzero$
and $j_{\ell_{1}}=j_{\ell_{2}}$, $\nn_{\ell_{1}}= \nn_{\ell_{2}}$,
or (ii) $j_{\ell_{1}}=3$, $\nn_{\ell_{1}}=\vzero$ and
$j_{\ell_{2}}=j_{\ell_{0}}$, $\nn_{\ell_{2}} = \nn_{\ell_{0}}\neq \vzero$
or (iii) $j_{\ell_{2}}=3$, $\nn_{\ell_{2}}=\vzero$ and
$j_{\ell_{1}}=j_{\ell_{0}}$, $\nn_{\ell_{1}} = \nn_{\ell_{0}}\neq \vzero$.
Moreover in case (i) one has to take the complex conjugate of all
propagators, node factors and endpoint factors of the subtree
with root line $\ell_{2}$. In all cases $k_{v}=0$ and $\nn_{v}=\vzero$,
so that the conservation law is obeyed also in this case;
cf. Figure \ref{fig:3.7}. The corresponding \textit{node factor} is
\begin{equation}
F_{v} := \begin{cases}
(1/2) (-1)^{j} \delta_{\nn_{v},\vzero} \delta^{\star}_{v} ,
& p_{v} =2 , \hbox{ case (i)} , \\
(1/2) (-1)^{j+1} i \delta_{\nn_{v},\vzero} \delta^{\star}_{v} ,
& p_{v}=2 , \hbox{ cases (ii) and (iii)} ,
\end{cases}
\label{eq:3.23} \end{equation}
where $\delta^{\star}_{v}$ recalls the constraints on the labels of
the entering and exiting lines of $v$, which are detailed above
and illustrated in Figure \ref{fig:3.8}.
The factor $1/2$ in the second line of (\ref{eq:3.23}) aims
to avoid overcountings of trees.

The nodes with branching number 1 can be only black points,
because there are no more white squares.
Hence (\ref{eq:3.15}) must be replaced with
\begin{equation}
F_{v} := f_{j_{\ell_{v}},j_{\ell_{v}'},\nn_{v}}
\left( 1 - \delta_{j_{\ell_{v}},3} \right) +
f_{1,j_{\ell_{v}'},\nn_{v}} \delta_{j_{\ell_{v}},3} , \qquad  p_{v}=1 ,
\label{eq:3.24} \end{equation}
which represents the \textit{node factor} of any node $v$ with $p_{v}=1$.
The corresponding order label is $k_{v}=1$, always.
A line $\ell$ exiting from a node $v$ can have
also $\nn_{\ell}=\vzero$ when $j_{\ell}=2$.

All endpoints $v$ have, by construction, $k_{v}=1$, and are
drawn as bullets coloured with black if $\nn_{\ell_{v}}\neq\vzero$
and coloured with white if $\nn_{\ell_{v}}=\vzero$,
in the latter case one must have $j_{\ell_{v}}=2$,
as $a_{\vzero}^{(1)}=0$; see (\ref{eq:3.11}).
The \textit{endpoint factor} of the endpoint $v$ is given by
\begin{equation}
F_{v} := f_{j_{\ell_{v}}1,\nn_{\ell_{v}}} ,
\label{eq:3.25} \end{equation}
which replaces (\ref{eq:3.16}). If $v$ is a white bullet
then necessarily $j_{\ell_{v}}=2$.

Finally, with the new rules, the propagator of any
line $\ell$ is given by
\begin{equation}
g_{\ell} := \begin{cases}
-i/\oo\cdot\nn_{\ell} ,
& \nn_{\ell}\neq \vzero , \; j_{\ell}=1 , \\
-i/(\oo\cdot\nn_{\ell}+2\la_{0}) ,
& \nn_{\ell}\neq \vzero , \; j_{\ell}=2 , \\
1 , & \nn_{\ell}= \vzero , \; j_{\ell}=1, \\
-i/2\la_{0} , & \nn_{\ell}= \vzero , \; j_{\ell}=2, \\
i , & \nn_{\ell}= \vzero , \; j_{\ell}=3 ,
\end{cases}
\label{eq:3.26} \end{equation}
which replaces both (\ref{eq:3.17}) and (\ref{eq:3.18}).

\begin{figure}[ht]
\centering
\vskip.3truecm
\ins{107pt}{-28pt}{$j,\nn$}
\ins{129pt}{-30pt}{$\vzero$}
\ins{150pt}{-17pt}{$j,\nn$}
\ins{138pt}{-37pt}{$3,\vzero$}
\ins{163.5pt}{-47pt}{$\nn_{1}$}
\ins{164pt}{-56pt}{$j_{1},-\nn_{1}$}
\ins{195pt}{-65pt}{$-\nn_{1}$}
\ins{211pt}{-75pt}{$1,\vzero$}
\ins{236pt}{-84pt}{$\vzero$}
\ins{252pt}{-74pt}{$j_{2},\nn_{2}$}
\ins{232pt}{-93pt}{$j_{2},\nn_{2},*$}
\ins{270pt}{-65.5pt}{$\nn_{3}$}
\ins{291pt}{-56pt}{$j_{4},\nn_{4}$}
\ins{311pt}{-48pt}{$\nn_{4}$}
\ins{330pt}{-38pt}{$2,\vzero$}
\includegraphics{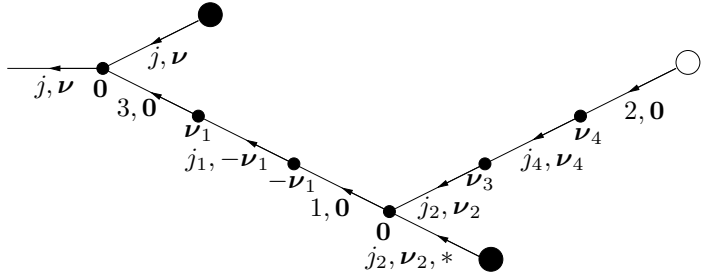}
\caption{An example of tree of order $k=7$ with $6$ nodes,
$3$ endpoints and $9$ lines.
All endpoints and all nodes with branching number $1$
have order $1$, while the nodes with branching number $2$
have order $0$ (hence it is useless to write the orders explicitly).
In principle the mode labels of the nodes with branching number $2$
and the labels of the lines coming out from the endpoints which are
white bullets could be omitted, as they are uniquely determined.
The conservation law for the momenta has been taken into
account explicitly, except for $\nn_{2}=\nn_{3}+\nn_{4}$.
One has $|V_{1}(\theta)|=4$
and $|V_{2}(\theta)|=|L_{0}(\theta)|=2$, so that $k_{\theta}=7$.}
\label{fig:3.9}
\end{figure}

An example of tree with the new rules is given in
Figure \ref{fig:3.9}. The order labels are not shown, for simplicity
(as the are identically $1$, except for the nodes with
branching number 2, which have order label $0$).

We still denote by $V(\theta)$ and $L(\theta)$ the number of nodes
and lines in $\theta$. Define also $E(\theta)$ the number of
endpoints of $\theta$, and set $P(\theta)=V(\theta)\cup E(\theta)$.
Furthermore call $V_{p}(\theta)$, $p=1,2$, the set of nodes
$v\in V(\theta)$ with branching number $p_{v}=p$, and
$L_{0}(\theta)$ the set of lines $\ell\in L(\theta)$
with $\nn_{\ell}=\vzero$ which do not come out from endpoints.
Then one has $|L_{0}(\theta)|=|V_{2}(\theta)|$.

We say that two trees are equivalent if they can be transformed
into each other by continuously deforming the lines in such a way
that the latter do not cross each other. Define $\Theta_{k,j,\nn}$
the set of inequivalent trees with labels $j,\nn$ associated with
the root line, and with $k_{\theta}=\sum_{v\in P(\theta)}k_{v}=
|P(\theta)|-|V_{2}(\theta)|=k$. The number of inequivalent
trees in $\Theta_{k,j,\nn}$ with fixed assignments
of the mode labels $\{\nn_{v}\}_{v\in P(\theta)}$
can be bounded by a constant to the power $k$: indeed 
a tree of order $k$ has $P(\theta) \le 2k$, so that the number
of unlabelled trees of order $k$ can be bounded by the number of
random walks with $4k$ steps, i.e. by $2^{4k}$, and all labels
except the mode labels assume a finite number of values.

We can summarise the considerations above into
the following formal statement.

\begin{lemma}\label{lem:8}
Let $\oo$ be a Bryuno vector. One has
\begin{eqnarray}
& & u^{(k)}_{j,\nn} = \sum_{\theta\in\Theta_{k,j,\nn}}
\Val (\theta) , \qquad k \ge 1, \; j=1,2, \nonumber \\
& & \mu^{(k)} = \sum_{\theta\in\Theta_{k,3,\vzero}} \Val (\theta) ,
\qquad k \ge 1 ,
\label{eq:3.27} \end{eqnarray}
with the \textit{tree value} $\Val(\theta)$ given by
\begin{equation}
\Val(\theta) = \Big( \prod_{\ell\in L(\theta)} g_{\ell} \Big)
\Big( \prod_{v \in P(\theta)} F_{v} \Big) .
\label{eq:3.28} \end{equation}
with the propagators $g_{\ell}$ given by (\ref{eq:3.26}),
and the factors $F_{v}$ given by (\ref{eq:3.23}),
(\ref{eq:3.24}) and (\ref{eq:3.25}). One has
$\mu^{(k)} \in i \RRR$ for all $k \ge 1$.
\end{lemma}

Even if (\ref{eq:3.28}) looks the same as (\ref{eq:3.20}),
the meaning of the symbols is different.

The formal series (\ref{eq:3.27}) is well defined, as it is easy
to check, but to order $k$, in general, we obtain for $\Val(\theta)$
bounds growing like $k!$ to some positive powers, so that summability
is prevented if we try to estimate the series (\ref{eq:3.27})
by taking the absolute values of the tree values.
To give a meaning to the formal series, we have to exploit
some remarkable cancellations between the tree values.
This can be showed by introducing a suitable resummation 
criterion of the series, which lead to a new series in which
to any order $k$ each tree value can be bounded
proportionally to a constant to the power $k$.
This will be done next.

\zerarcounters
\section{Renormalised series} 
\label{sec:4}

Consider a tree $\theta$, and suppose that each line $\ell$
carries a further label $n_{\ell} \in \ZZZ_{+} \cup \{-1\}$,
the \textit{scale} label. We say that a connected set of lines
and nodes $T \subset L(\theta)$ is a \textit{cluster} on scale $n_{T}$
if (i) all lines in $T$ have scales no smaller than $n_{T}$,
(ii) at least one line in $T$ is on scale $n_{T}$, and
(iii) it is maximal (which means that the lines connected to $T$
but not belonging to it are on scales less than $n_{T}$).
If $T$ contains only one node (and no lines) we set $n_{T}=-1$,
as in the case in which all lines in $T$ are on scale $-1$.

If $\theta$ is a linear tree then all clusters
have only one entering line, while in nonlinear trees clusters
can have any number of entering lines. On the contrary a cluster,
in both linear and nonlinear trees, can have only either zero or
one exiting line. We call \textit{external lines} of a cluster $T$
the lines which are either entering or exiting lines for $T$.

We say that the cluster $T$ is a \textit{self-energy cluster} if
(i) $T$ has only one entering line and only one exiting line,
(ii) the entering line carries the same momentum and
component label as the exiting line, and
(iii) no line along the path of lines connecting the external lines
has vanishing momentum.

\begin{figure}[ht]
\centering
\vskip.3truecm
\ins{84pt}{-81pt}{$j,\nn$}
\ins{111pt}{-84pt}{$\vzero$}
\ins{130pt}{-73pt}{$3,\vzero$}
\ins{146pt}{-66pt}{$\nn_{0}$}
\ins{168pt}{-53pt}{$j_{0},-\nn_{0}$}
\ins{120pt}{-91pt}{$j,\nn$}
\ins{145pt}{-102pt}{$\nn_{1}$}
\ins{152.5pt}{-112pt}{$j_{1},\!\!-\nn_{1}$}
\ins{176pt}{-119pt}{$-\nn_{1}$}
\ins{195pt}{-128pt}{$j,\nn$}
\ins{219pt}{-138pt}{$\nn$}
\ins{226pt}{-146pt}{$1,\vzero$}
\ins{255.5pt}{-156pt}{$\vzero$}
\ins{270pt}{-145pt}{$j_{2},\nn_{2}$}
\ins{252pt}{-166pt}{$j_{2},\nn_{2},*$}
\ins{291pt}{-138pt}{$\nn_{3}$}
\ins{306pt}{-130pt}{$j_{4},\nn_{4}$}
\ins{328pt}{-120pt}{$\nn_{4}$}
\ins{348pt}{-109pt}{$2,\vzero$}
\ins{168pt}{-91pt}{$T$}
\includegraphics{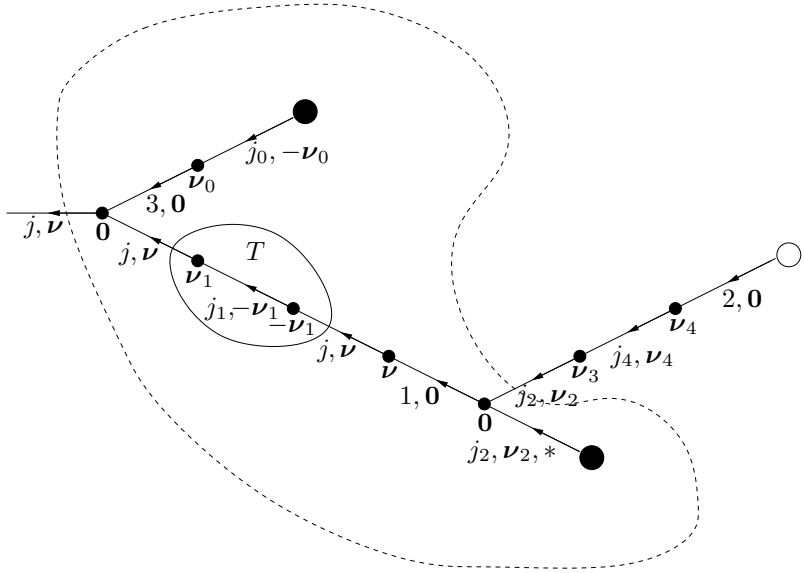}
\caption{Example of self-energy cluster. Let $T$ be the set
of nodes and lines inside the solid line, i.e. the set
consisting of the line $\ell$ with momentum $-\nn_{1}$
and of the two nodes $v_{1}$ and $v_{2}$, with mode labels
$\nn_{1}$ and $-\nn_{1}$, respectively, connected by such a line.
Then $T$ is a self-energy cluster if the scale of the line $\ell$
is strictly less than the scales of both the line $\ell_{2}$
entering $v_{2}$ and the line $\ell_{1}$ exiting $v_{1}$, i.e. if
$n_{\ell} < \min\{n_{\ell_{1}},n_{\ell_{2}}\}$. In such
a case $\ell_{1}$ and $\ell_{2}$ become the external lines of $T$.
The set of nodes and lines inside the dotted line cannot
be a self-energy cluster, even if it is a cluster and $\nn=\nn_{2}$,
because the path of lines between the external lines
contains a line with vanishing momentum.}
\label{fig:4.1}
\end{figure}

A self-energy cluster by construction can contain
other self-energy clusters. We say that a self-energy cluster
is a \textit{renormalised self-energy cluster}
if it does not contain any other self-energy clusters.
We say that a tree $\theta$ is a \textit{renormalised tree}
if it does not contain any self-energy clusters.
Given a self-energy cluster $T$, denote by $V(T)$, $E(T)$ and $L(T)$
the set of nodes, the set of endpoints and the set of lines,
respectively, contained in $T$, and set $P(T)=V(T)\cup E(T)$.
Call $V_{p}(T)$ the set of nodes $v \in V(T)$ with $p_{v}=p$,
and $L_{0}(T)$ the set of lines $\ell\in L(T)$ with $\nn_{\ell}=\vzero$
which do not come out from endpoints. Set $k_{T}=|P(T)|-|V_{2}(T)|$. 
An example of self-energy cluster is given in Figure \ref{fig:4.1}.

Define the \textit{self-energy value} $\VV_{T}(\oo\cdot\nn)$ as
\begin{equation}
\VV_{T}(\oo\cdot\nn) = \eps^{k_{T}}
\Big( \prod_{\ell\in L(T)} g_{\ell}^{\RR} \Big)
\Big( \prod_{v \in P(T)} F_{v} \Big) ,
\label{eq:4.1} \end{equation}
with the factors $F_{v}$ defined as in Section \ref{sec:3}
and the \textit{renormalised propagators} $g^{\RR}_{\ell}$
still to be defined.

The renormalised self-energy clusters can be of two kinds:
those in which both external lines are attached to the same node,
and those in which there is a nontrivial path of lines
connecting the external lines. Those of the first
type can be seen as obtained from the expansion
of the white square representing a node of in a linear tree.

\begin{figure}[ht]
\centering
\vskip.3truecm
\ins{139pt}{-12pt}{$T$}
\ins{152pt}{-14pt}{$=$}
\ins{137pt}{-84pt}{$T'$}
\ins{152pt}{-86pt}{$=$}
\ins{135pt}{-156pt}{$T''$}
\ins{152pt}{-158pt}{$=$}
\ins{173pt}{-18pt}{$\ell_{\rm out}$}
\ins{270pt}{-18pt}{$\ell_{\rm in}$}
\ins{173pt}{-90pt}{$\ell_{\rm out}$}
\ins{222pt}{-101pt}{$\ell_{\rm in}$}
\ins{173pt}{-163pt}{$\ell_{\rm out}$}
\ins{221pt}{-135pt}{$\ell_{\rm in}$}
\ins{208pt}{-19pt}{$v_{\rm out}$}
\ins{234pt}{-19pt}{$v_{\rm in}$}
\ins{206pt}{-77pt}{$v_{0}$}
\ins{260pt}{-46pt}{$v_{\rm in}$}
\ins{239pt}{-75pt}{$v_{\rm out}$}
\ins{226pt}{-81pt}{$\ell_{0}$}
\ins{206pt}{-165pt}{$v_{0}$}
\ins{263pt}{-197pt}{$v_{\rm in}$}
\ins{236pt}{-183pt}{$v_{\rm out}$}
\ins{222pt}{-175pt}{$\ell_{0}$}
\includegraphics{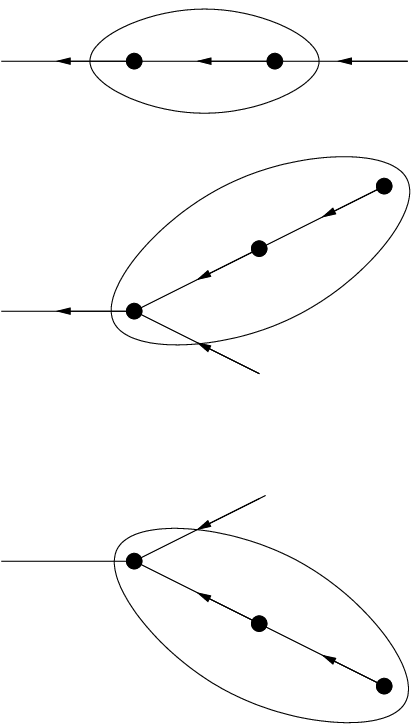}
\caption{Examples of renormalised self-energy clusters belonging
to the same equivalence class: $T$ is a renormalised
self-energy cluster of the second type, while $T'$ and $T''$ are 
renormalised self-energy clusters of the first kind.
The self-energy values of $T'$ and $T''$ are equal to each other:
in fact the trees containing such renormalised self-energy clusters
can be obtained from each other by
permuting the entering lines of $v_{0}$. The external lines
$\ell_{\rm in}$ and $\ell_{\rm out}$ do not belong to the
self-energy clusters, and have been drawn only to help
visualising the structure of the self-energy clusters.}
\label{fig:4.2}
\end{figure}

Consider a renormalised self-energy cluster $T$ of the second kind.
Call $v_{\rm in}$ and $v_{\rm out}$ the nodes which the entering line
$\ell_{\rm in}$ and the exiting line $\ell_{\rm out}$ of $T$,
respectively, are attached to. Then add a further node $v_{0}$
and a further line $\ell_{0}$ and consider the set $\tilde T$, with
$V(\tilde T)=V(T)\cup \{v_{0}\}$ and $L(\tilde T)=L(T)\cup\{\ell_{0}\}$,
constructed as follows. Detach the line $\ell_{\rm out}$
from $v_{\rm out}$ add attach it to the node $v_{0}$, and
connect the node $v_{0}$ to the node $v_{\rm out}$ through the
line $\ell_{0}$ (oriented from $v_{\rm out}$ to $v_{0}$).
Finally detach the line $\ell_{\rm in}$ from $v_{\rm in}$ and
reattach it to the node $v_{0}$ (so that $p_{v_{0}}=2$).
The last operation can be performed in two ways ($\ell_{\rm in}$
can be above or below $\ell_{0}$), hence it generates two
renormalised self-energy clusters $T'$ and $T''$. We 
call, shortly, \textit{shift operation} the mechanism
described above; cf. Figure \ref{fig:4.2}.

\begin{lemma}\label{lem:9}
For each renormalised self-energy cluster $T$ of the second kind
there is one and only one pair of renormalised self-energy clusters
$T'$ and $T''$ of the first kind which can be obtained from $T$
through the shift operation.
\end{lemma}

\prova The proof is a simple application of the diagrammatic rules
described in Section \ref{sec:3}.$\EP$

This allows us to introduce a notion of equivalence between renormalised
self-energy clusters. Then we can consider the renormalised
self-energy clusters as triples of equivalent
renormalised self-energy clusters $\{T,T',T''\}$.

Assume that $\oo$ be a Bryuno vector. Define
\begin{equation}
C_{0} = \sum_{n=0}^{\infty} 2^{n(d-1)} \al_{n} , \qquad
\al_{n} = \al_{n}(\om) := \inf_{|\nn| \le 2^{n}} |\oo\cdot\nn| ,
\label{eq:4.2} \end{equation}
and set $\al_{n}=C_{0}\g_{n}$. If the sum in (\ref{eq:4.2}) diverges,
redefine $C_{0}$ by writing $2^{n(d-2)}$ instead of $2^{n(d-1)}$
(so that convergence is assured because $\al_{n} \le |\oo| 2^{-n(d-1)}$,
by Dirichlet's theorem \cite{S}), and replace $\g_{n}$
with $\g_{n}2^{-n}$ in the following multiscale decomposition
-- see the definition of the compact support functions $\chi_{n}$
after (\ref{eq:4.4}), -- and in the
Diophantine conditions (\ref{eq:4.14}).

Note that $n'>n$ implies $\g_{n'} \le \g_{n}$,
while $\g_{n'}<\g_{n}$ implies $n'>n$.

Set $\ZZZ^{d}_{*}=\ZZZ^{d}\setminus\{\vzero\}$, and define
\begin{equation}
n(\nn) = \{ n \in \ZZZ_{+} : 2^{n-1} < |\nn| \le 2^{n} \}
\label{eq:4.3} \end{equation}
for all $\nn\in\ZZZ^{d}_{*}$.

Let $\psi(x)$ a non-decreasing $C^{\infty}(\RRR)$ function
defined in $\RRR$, such that
\begin{equation}
\psi(x) = \left\{
\begin{array}{ll}
1 \, , & \text{for } |x| \geq C_{1} \, , \\
0 \, , & \text{for } |x| \leq C_{1} /2 \, ,
\end{array} \right.
\label{eq:4.4} \end{equation}
with the constant $C_{1} \le C_{0}$ to be defined later.
Set also $\chi(x) := 1-\psi(x)$, and define, for all $n \in \ZZZ_{+}$,
$\chi_{n}(x) := \chi(\beta^{-1}\g_{n}^{-1} x)$ and
$\psi_{n}(x) := \psi(\beta^{-1}\g_{n}^{-1}x)$, with $\beta=1/4$.

Define
\begin{equation}
\Delta_{0}(x) = \left( \frac{1}{2} 
\left( \frac{1}{x^{2}} + \frac{1}{(x+2\la_{0})^{2}} \right)
\right)^{-1/2} ,
\label{eq:4.5} \end{equation}
and, setting $\MM^{[0]}_{1}(x) := 0$ and
$\MM^{[0]}_{2}(x) := \la_{0}$, define for $n \ge 1$ and $j=1,2$
\begin{eqnarray}
\MM^{[\le n]}_{j}(x)
& = & \sum_{p=0}^{n} \MM^{[p]}_{j}(x) , \nonumber \\
\MM^{[n]}_{j}(x) & = & \chi_{0}(\Delta_{0}(x)) \ldots
\chi_{n-1}(\Delta_{0}(x)) M^{[n]}_{j}(x) , \nonumber \\
M^{[n]}_{j}(x) & = & \frac{i}{2}\sum_{k=1}^{\infty}
\sum_{T \in \SSS_{k,j,n-1}} \VV_{T}(x) ,
\label{eq:4.6} \end{eqnarray}
where $\SSS_{k,j,n}$ is the set of all renormalised self-energy clusters
$T$ on scale $n$ with $|P(T)|-|V_{2}(T)|=k$
and with component label $j$ associated with both external lines.
For $n=0$ we interpret $\MM^{[\le 0]}_{j}(x)=\MM^{[0]}_{j}(x)$.
One has
\begin{equation}
\min\{|x|,|x+2\la_{0}|\}
\le \Delta_{0}(x) \le \sqrt{2} \min\{|x|,|x+2\la_{0}|\} .
\label{eq:4.7} \end{equation}

Then the \textit{renormalised propagator} is defined as
$g^{\RR}_{\ell}=g_{\ell}$ if $\nn_{\ell}=\vzero$
and $g^{\RR}_{\ell}=g^{[n_{\ell}]}_{j_{\ell}}
(\oo\cdot\nn_{\ell})$ if $\nn_{\ell} \neq \vzero$, with
\begin{equation}
g^{[n]}_{j}(x) = - i
\frac{\chi_{0}(\Delta_{0}(x)) \ldots \chi_{n-1}(\Delta_{0}(x))
\psi_{n}(\Delta_{0}(x))}{x+2 \MM^{[\le n]}_{j}(x)} ,
\label{eq:4.8} \end{equation}
so that we see that $g^{[n]}(x) \neq 0$ implies
\begin{eqnarray}
\frac{1}{2} \beta \g_{n} C_{1} \le
\Delta_{0}(x) \le \beta \g_{n-1} C_{1} .
\label{eq:4.9} \end{eqnarray}
We associate a scale label $n_{\ell}$ also with lines with
vanishing momentum, by setting $n_{\ell}=-1$. 

Note that $\MM^{[\le n]}_{j}(x)$  is defined in terms
of propagators on scales $n'<n$, hence in terms of
$\MM^{[n']}_{j'}(x')$, with $n'<n$: this means that (\ref{eq:4.6})
provides a recursive definition of $\MM^{[\le n]}_{j}(x)$,
hence it makes sense. Note also that self-energy clusters
on scale $-1$ (in particular those consisting of a single node)
are not taken into account in (\ref{eq:4.6});
this will be motivated by Lemma \ref{lem:10} below.

Define the \textit{tree value} $\Val(\theta)$ as
\begin{equation}
\Val(\theta) = \Big( \prod_{\ell\in L(\theta)} g^{\RR}_{\ell} \Big)
\Big( \prod_{v \in P(\theta)} F_{v} \Big) .
\label{eq:4.10} \end{equation}
Then, if $\Theta^{\RR}_{j,k,\nn}$ is the set of inequivalent
renormalised trees with labels $j,\nn$ associated with the root line
and with $|P(\theta)|-|V_{2}(\theta)|=k$, set
\begin{equation}
u_{j,\nn}^{[k]} = \sum_{\theta\in\Theta^{\RR}_{j,k,\nn}}
\Val(\theta) ,
\label{eq:4.11} \end{equation}
with $u_{3,\vzero}^{[k]}:=\mu^{[k]}$,
and define the function $\overline u(t)=(u_{1}(t),u_{2}(t))$ as 
\begin{equation}
\overline{u}_{j}(t) = \sum_{k=1}^{\infty} \eps^{k}
u^{[k]}_{j}(t) , \qquad
u^{[k]}_{j}(t) = \sum_{\nn\in\ZZZ^{d}} {\rm e}^{i\nn\cdot\oo t}
u_{j,\nn}^{[k]} ,
\label{eq:4.12} \end{equation}
and the \textit{counterterm} $\overline \mu$ as
\begin{equation}
\overline \mu = \sum_{k=1}^{\infty} \eps^{k} \mu^{[k]} ,
\label{eq:4.13} \end{equation}
that we call the \textit{renormalised series}
for $u(t)$ and $\mu$, respectively.

\begin{lemma}\label{lem:10}
The self-energy clusters on scale $-1$ have values which cancel out
exactly when summed together, hence there is no contributions
arising from them to $\MM^{[\le n]}_{j}(x)$.
\end{lemma}

\begin{figure}[ht]
\centering
\ins{129pt}{-36pt}{$j,\nn$}
\ins{152pt}{-38pt}{$\vzero$}
\ins{172pt}{-36pt}{$j,\nn$}
\ins{227pt}{-30.5pt}{$+$}
\ins{272pt}{-36pt}{$j,\nn$}
\ins{296.5pt}{-38pt}{$\vzero$}
\ins{312pt}{-36pt}{$j,\nn$}
\ins{319pt}{-17pt}{$3,\vzero$}
\includegraphics{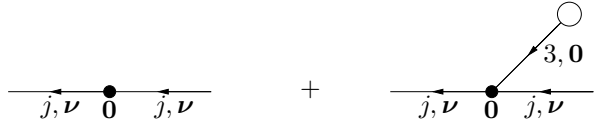}
\vskip.5truecm
\caption{Self-energy clusters on scale $-1$
contributing to $\MM^{[\le n]}_{j}(\oo\cdot\nn)$.
The external lines do not belong to the
self-energy clusters, and have been drawn only to help
visualising the structure of the self-energy clusters.}
\label{fig:4.3}
\end{figure}

\prova The self-energy clusters on scale $-1$ are those
represented in Figure \ref{fig:4.3}. Hence they would contribute
to $\MM^{[\le n]}_{j}(x)$ a value
$f_{11,\vzero} + i\mu^{[1]}$ for $j=1$ and
$f_{22,\vzero} - i\mu^{[1]}$ for $j=2$.
By the very definition of $\mu^{[1]}$ one has
$i\mu^{[1]}=-f_{11,\vzero}$, so that
$f_{11,\vzero} + i\mu^{[1]}=0$ for $j=1$. For $j=2$ one has
$f_{22,\vzero} - i\mu^{[1]}=f_{22,\vzero}+f_{11,\vzero}=0$,
where we used that $f\in\gotm$, so that $\tr f=0$.$\EP$

For higher values of $n$, $\MM^{[n]}_{1}(x)$ and $\MM^{[n]}_{2}(x)$
are no longer equal to each other. However, we shall see that there is
a deep symmetry yielding $\MM^{[n]}_{1}(0)=-\MM^{[n]}_{2}(-2\la_{0})$
(cf. Lemma \ref{lem:15}). Moreover, the cancellation mechanism which
leads to Lemma \ref{lem:10} still works for any $n$, and implies
partial cancellations, as Lemma \ref{lem:16} will show.

Assume the Diophantine conditions
\begin{eqnarray}
\left| \oo\cdot \nn \right| > C_{1} \g_{n(\nn)} , \nonumber \\
\left| \oo\cdot \nn \pm 2 \la_{0} \right| >
C_{1} \g_{n(\nn)} ,
\label{eq:4.14} \end{eqnarray}
for all $\nn\in\ZZZ^{d}_{*}$ and all $n\ge 0$. For $C_{1} \le C_{0}$ the
conditions in the first line are automatically satisfied by definition.
The condition in the second line, called the \textit{(first)
Melnikov condition}, instead, have to be explicitly required with the
constant $C_{1}$ -- the same as in (\ref{eq:4.4}) -- still to be fixed.

Let $\Lambda_{0}$ be the set in which $\la_{0}$ varies,
and call $\Lambda_{0}^{*}$ the subset of values $\la_{0}\in\Lambda_{0}$
for which the conditions (\ref{eq:4.14}) are satisfied.
Of course $\Lambda_{0}$ has to be such that for $\la_{0}\in\Lambda_{0}$
one has $\la=\la_{0}+\mu \in [a,b]$, but for the time being
we ignore such a constraint.

\zerarcounters
\section{Convergence of the renormalised series}
\label{sec:5}

In this Section we assume that $\la_{0}\in\Lambda_{0}^{*}$. Hence
the Diophantine conditions (\ref{eq:4.14}) are satisfied. We want
to study the renormalised series for $u$ and $\mu$, with the aim
of showing first that they converge, so that the functions
$\overline u$ and $\overline \mu$ are well defined,
second that $\overline u$ solves the equations (\ref{eq:2.19})
provided one fixes $\mu=\overline\mu$ and both $\overline u$
and $\overline \mu$ are analytic in $\eps$, third that
the relative measure of the set $\Lambda_{0}^{*}$ is large.
Finally we have to check that
the last property implies that the set of values
of $\la$ in $[a,b]$ for which (\ref{eq:2.10}) is reducible also is of
large measure; this will be done in Section \ref{sec:6}.

We note since now that for any renormalised self-energy cluster $T$
one has $|L(T)|=|P(T)|-1$, so that
$|L(T)|-|L_{0}(T)|=k_{T}-1$. Moreover if $T\in \SSS_{k,j,n}$,
with $n \ge 0$, then $k_{T}\ge 2$, because there must be
at least one line on scale $n$.

In the following by saying that some property holds ``for $\eps$
small enough'' we mean that there exists a
constant $\eps_{0}$ (not necessarily the same
in all the statements) such that (i) $\eps_{0}C_{1}^{-1} \ll 1$,
and (ii) for $|\eps|<\eps_{0}$
that property is satisfied. Define also
\begin{equation}
|x+2\rho_{0}(x)| := \min\{ |x|,|x+2\la_{0}| \} ,
\label{eq:5.1} \end{equation}
so that $\rho_{0}(x)$ is either $0$ (if $x + \la_{0} \ge 0$)
or $\la_{0}$ (if $x+\la_{0} < 0$).

An important remark is that in the forthcoming Lemmata \ref{lem:11}
and \ref{lem:12} the results hold unchanged if, in (\ref{eq:4.9}),
we replace $\beta$ with $2\beta$ in the upper bound and $\beta$
with $\beta/2$ in the lower bound. Why this is important will be
explained in the proof of Lemma \ref{lem:16}.

\begin{lemma}\label{lem:11}
Let $\theta$ be a renormalised tree. Call $N_{n}(\theta)$
the number of lines in $L(\theta)$ on scale $n$. One has
\begin{equation}
N_{n}(\theta) \le K 2^{-n} M(\theta) , \qquad
M(\theta) := \sum_{v\in P(\theta)} |\nn_{v}|
\label{eq:5.2} \end{equation}
for a suitable constant $K$.
\end{lemma}

\prova First note that if $N_{n}(\theta) \neq 0$ then, by (\ref{eq:4.9}),
there exists a line $\ell\in L(\theta)$ such that $C_{1}\g_{n-1} >
\beta C_{1} \g_{n-1} \ge |x+2\rho_{0}(x)| >
C_{1} \g_{n(\nn_{\ell})}$, hence $n(\nn_{\ell})\ge n$,
thence $M(\theta) \ge |\nn_{\ell}| > 2^{n-1}$.

Then we prove by induction that
\begin{equation}
N_{n}(\theta) \neq 0 \quad \Longrightarrow \quad
N_{n}(\theta) \le 2 2^{-n} M(\theta) - 1 .
\label{eq:5.3} \end{equation}
If the root line of $\theta$ is not on scale $n$ the bound
(\ref{eq:5.3}) follows by induction. If the root line $\ell$ of $\theta$
is on scale $n$ consider the lines
$\ell_{1},\ldots,\ell_{p}$ on scales $\ge n$ such that
no line along the paths connecting any of them to the root line
is on scale $\ge n$. If $p \ge 2$ again the bound follows inductively.
If $p=1$ call $\theta_{1}$ the subtree with root line $\ell_{1}$,
and call $T$ the set of points and lines between $\ell_{1}$
and $\ell$ (that is which precede $\ell$ but not $\ell_{1}$).
Denote by $P(T)$ the set of points in $T$, and
define $M(T):=\sum_{v\in P(T)}|\nn_{v}|$.
Call $\nn$ and $\nn'$ the momenta associated with
$\ell$ and $\ell_{1}$, respectively, and set $x=\oo\cdot\nn$ and
$x'=\oo\cdot\nn'$. One has $N_{n}(\theta)=1+N_{n}(\theta_{1})$,
and both $|x + 2 \rho_{0}(x)|$ and $|x' + 2 \rho_{0}(x')|$
are less than $\beta C_{1}\g_{n-1}$, so that
$|(x-x') + 2 (\rho_{0}(x)-\rho_{0}(x'))| \le
2 \beta C_{1} \g_{n-1} < C_{1} \g_{n-1}$.

If there is (at least) a line $\ell'$ with $\nn_{\ell'}=\vzero$ along
the path of lines between the external lines $\ell$ and $\ell_{1}$,
then there exist two disjoint sets $T_{1}$ and $T_{2}$, with $P(T)=
P(T_{1})\cup P(T_{2})$ and $L(T)=L(T_{1})\cup L(T_{2}) \cup \{\ell'\}$,
such that both $M(T_{1})$ and $M(T_{2})$ are greater than $|\nn|$.
Since $\ell$ is on scale $n$ one has $|\nn|>2^{n-1}$,
so that $M(T) \ge \max\{M(T_{1}),M(T_{2})\} \ge 2^{n-1}$.
If there is no line with zero momentum between the
external lines, then $\nn\neq\nn'$, otherwise $T$
would be a renormalised self-energy cluster. Therefore by the
second Diophantine conditions (\ref{eq:4.14}), one obtains
$n(\nn-\nn') \ge n$, so that $M(T) \ge 2^{n-1}$ also in such a case.

Hence, by the inductive hypothesis
\begin{equation}
N_{n}(\theta) \le 1 + \left( 2 2^{-n} M(\theta_{1}) - 1 \right)
\le 1 - 2 2^{-n} M(T) + \left( 2 2^{-n} M(\theta) - 1 \right) \le
2 2^{-n} M(\theta) -1 ,
\label{eq:5.4} \end{equation}
and the bound (\ref{eq:5.3}) follows.$\EP$

\begin{lemma}\label{lem:12}
Let $T$ be a renormalised self-energy cluster. Call $N_{n}(T)$ the
number of lines in $L(T)$ on scale $n$, with $n \le n_{T}$. One has
\begin{equation}
N_{n}(T) \le K 2^{-n} M(T) , \qquad
M(T) := \sum_{v\in P(T)} |\nn_{v}| > 2^{n_{T}-1} ,
\label{eq:5.5} \end{equation}
with the same constant $K$ as in (\ref{eq:5.2}).
\end{lemma}

\prova We first prove the bound on $M(T)$. By construction
$T$ must contain at least a line $\ell$ on scale $n_{T}$, so that
$|x_{\ell}+2\rho_{0}(x_{\ell})| \le \beta C_{1} \g_{n_{T}-1}$,
with $x_{\ell}=\oo\cdot\nn_{\ell}$.
Write $\nn_{\ell}=\nn_{\ell_{0}}^{0}+\sigma_{\ell}\nn$,
where $\nn$ is the momentum associated with the entering line
of $T$ and $\sigma_{\ell}=0,1$, and set $x=\oo\cdot\nn$
and $x_{\ell}^{0}=\oo\cdot\nn_{\ell}^{0}$. The
entering line of $T$ has scale strictly larger than $n_{T}$,
so that $|x+2\rho_{0}(x)| \le \beta C_{1} \g_{n_{T}-1}$.
If $M(T) \le 2^{n_{T}-1}$ then $|\nn_{\ell}^{0}|\le M(T) \le 2^{n_{T}-1}$,
hence $n(\nn_{\ell}^{0}) \le n_{T}-1$, so that $|x_{\ell}^{0}+2
\rho_{0}(x_{\ell}^{0})| > C_{1} \g_{n(\nn_{\ell}^{0})}
\ge C_{1} \g_{n_{T}-1}$, by the Diophantine conditions (\ref{eq:4.14}).
Then one has
\begin{eqnarray}
C_{1} \g_{n_{T}-1} & > & |x_{\ell} + 2\rho_{0}(x_{\ell})|
+ \sigma_{\ell} |x + 2\rho_{0}(x)| \nonumber \\
& \ge & |x_{\ell}^{0} + 2 ( \rho_{0}(x_{\ell}) - \sigma_{\ell}
\rho_{0}(x) | > C_{1} \g_{n_{T}-1} ,
\label{eq:5.6} \end{eqnarray}
which leads to a contradiction.

Next we pass to the bound on $N_{n}(T)$.
Consider a subset $G_{0}$ of the lines of a tree $\theta$
between two lines $\ell_{\rm out}$ and $\ell_{\rm in}$
Set $G=G_{0}\cup \{\ell_{\rm in}\}\cup\{\ell_{\rm out}\}$.
Let $[n_{\rm in}],[n_{\rm out}]$ be the scales of the lines
$\ell_{\rm out}$ and $\ell_{\rm in}$, respectively,
and suppose that $n_{\rm in},n_{\rm out}\ge n$, while all lines in
$G_{0}$ (if any) have scales $n' \le n_{T}-1$. Note that in general
$G_0$ is not even a cluster unless $n_{\rm in},n_{\rm out}\ge n_{T}$.
Then we can prove that if $N_{n}(G_{0})\neq0$ then
$N_{n}(G_{0}) \le 2 2^{-n} \sum_{v \in P(G_0)}|\nn_{v}|-1$,
where $P(G_{0})$ is the set of points preceding $\ell_{\rm out}$
and following $\ell_{\rm in}$.

If $G_{0}$ has no lines then the mode $\nn_{0}$ of the (only)
node between $\ell_{\rm out}$ and $\ell_{\rm in}$
is such that $|\nn_{0}| \ge 2^{n-1}$,
by the second Diophantine conditions (\ref{eq:4.14}),
and the statement is true. Hence we proceed inductively
on the number of lines in $G_{0}$.
If no line of $G_{0}$ on the path ${\cal P}(G)$ connecting the
external lines of $G$ has scale $n$ then the lines in $G_{0}$
on scale $n$ (if any) belong to trees with root on
${\cal P}(G)$, and the statement follows from the
bound (\ref{eq:5.3}) for trees given in the proof of Lemma \ref{lem:11}.
If there is a line $\ell\in {\cal P}(G)$ on scale $n$,
then call $G_{1}$ and $G_{2}$ the disjoint subsets of $G$ such that
$G_{1}\cup G_{2} \cup \{\ell\} = G$. Then $G_1\cup\{\ell\}$ and
$G_2\cup\{\ell\}$ have the same structure of $G$ itself, but each has
less lines. Hence, again the inductive assumption yields the result.

Therefore, as a particular case, by choosing $G_0=T$, with
$T\in\SSS_{k,j,n_{T}-1}$, the bound for $N_{n}(G_{0})$
implies the bound on $N_{n}(T)$ we are looking for.$\EP$

\begin{lemma}\label{lem:13}
Assume that the propagators $g^{[p]}_{j}(x)$
can be uniformly bounded for all $0\le p \le n-1$ as
\begin{equation}
\left| g^{[p]}_{j}(x) \right| \le K_{1} C_{1}^{-1} \g_{p}^{-1} ,
\label{eq:5.7} \end{equation}
for some positive constant $K_{1}$. Then one has
\begin{equation}
\left| \VV_{T}(\oo\cdot\nn) \right| \le
|\eps|^{k_{T}} D_{1}^{k_{T}} C_{1}^{- (k_{T}-1)}
\g_{m_{0}}^{-k_{T}}
{\rm e}^{-\ka_{0} M(T)/2} ,
\label{eq:5.8} \end{equation}
for a suitable constant $D_{1}$. If also
the derivatives of the propagators are bounded uniformly as
\begin{equation}
\left| \partial_{x} g^{[p]}_{j}(x) \right|
\le K_{2} C_{1}^{-2} \g_{p}^{-3} ,
\label{eq:5.9} \end{equation}
for some positive constant $K_{2}$, one has also
\begin{equation}
\left| \frac{{\rm d}}{{\rm d}x} \left. \VV_{T}(x)
\right|_{x=\oo\cdot\nn} \right| \le
|\eps|^{k_{T}} D_{2}^{k_{T}} C_{1}^{-k_{T}}
\g_{m_{0}}^{-k_{T}-2} {\rm e}^{-\ka_{0}M(T)/2} ,
\label{eq:5.10} \end{equation}
for a suitable constant $D_{2}$.
\end{lemma}

\prova For any renormalised self-energy cluster $T$ consider the 
corresponding self-energy value (\ref{eq:4.1}). The product of factors 
$F_{v}$ can be bounded as
\begin{equation}
\prod_{v \in P(T)} F_{v} \le F_{0}^{k_{T}}
\prod_{v \in V_{1}(T) \cup E(T)} {\rm e}^{-\ka_{0}|\nn_{v}|} ,
\label{eq:5.11} \end{equation}
while the product of propagators can be bounded,
for any $m_{0}\in \NNN$, as
\begin{equation}
\prod_{v \in L(T)} g^{\RR}_{\ell} \le C_{1}^{- (k_{T}-1)}
\g_{m_{0}}^{-k_{T}} \exp \left( K
\sum_{n=m_{0}+1}^{\infty} \frac{1}{2^{n}} \log \frac{1}{\g_{n}}
M(T) \right) ,
\label{eq:5.12} \end{equation}
where the first bound (\ref{eq:5.5}) of Lemma \ref{lem:12}
has been used. If we choose $m_{0}$ such that
\begin{equation}
K \sum_{n=m_{0}+1}^{\infty} \frac{1}{2^{n}} \log \frac{1}{\g_{n}}
\le \frac{\ka_{0}}{12} ,
\label{eq:5.13} \end{equation}
then we obtain (\ref{eq:5.8}). Such $m_{0}$ exists
because $\oo$ is a Bryuno vector; cf. (\ref{eq:3.2}).

Call $\PP(T)$ the path of lines $\ell\in L(T)$ which are between
the external lines of $T$.
Then the derivative of $\VV_{T}(x)$ can be written as
\begin{equation}
\partial_{x} \VV_{T}(x) = \eps^{k_{T}}
\sum_{\ell\in \PP(T)} \partial_{x} g_{\ell}^{\RR}
\Big( \prod_{\ell'\in L(T) \setminus \ell} g_{\ell'}^{\RR} \Big)
\Big( \prod_{v \in P(T)} F_{v} \Big) ,
\label{eq:5.14} \end{equation}
so that, by reasoning as in the previous case,
using the bounds (\ref{eq:5.9}) and choosing again
$m_{0}$ as in (\ref{eq:5.13}), we obtain (\ref{eq:5.10}).$\EP$

\begin{lemma}\label{lem:14}
$\MM^{[\le n]}_{j}|\RRR$ is real for all $n\ge 0$ and $j=1,2$.
\end{lemma}

\prova The proof is by induction on $n$.
For $n=0$ the assertion is trivially satisfied.
Then assume that it holds for all $n'<n$.

Let $T$ be a renormalised self-energy cluster contributing to
$\MM^{[n]}_{j}(x)$ through (\ref{eq:4.6}). Denote by $v_{\rm in}$
and $v_{\rm out}$ the nodes in $V(T)$ which  the entering line
$\ell_{\rm in}$ and the exiting line $\ell_{\rm out}$ 
of $T$ are attached to, respectively. Call $\PP(T)$ the set
of lines and nodes between the external lines of $T$.

Together with $T$ consider also the renormalised self-energy cluster $T'$
obtained as follows. Detach the line $\ell_{\rm in}$ from $v_{\rm in}$
and attach it to the node $v_{\rm out}$, and
detach the line $\ell_{\rm out}$ from $v_{\rm out}$
and attach it to the node $v_{\rm in}$. Consistently, orient
all lines along the path $\PP(T)$ between the external lines of $T$
in the opposite direction, i.e. from $v_{\rm out}$ to $v_{\rm in}$.
Finally change the mode labels of all nodes along $\PP(T)$, i.e.
of all nodes $v\in V(\PP(T))$, if $V(\PP(T))$ denotes the set of nodes
along $\PP(T)$. The latter operation is possible because of the
following reason. Each line entering a node $v\in V(\PP(T))$ has
zero momentum: indeed for each node $v$ with branching number $p_{v}=2$
one of the three lines connected with $v$ must have zero momentum
(cf. Figure \ref{fig:3.8}), and by definition of self-energy cluster
such a line cannot lay on $\PP(T)$. Hence $\sum_{v \in V(\PP(T))}
\nn_{v}=\vzero$. Note also that each line entering a node
$v\in V(\PP(T))$ is the root line of a tree contributing to
$\mu^{[k_{v}]}$, for some $k_{v}$. Along the path $\PP(T)$
the propagators have not changed because of the operation above
(cf. the analogous discussion in the proof of Lemma \ref{lem:7}),
by the inductive hypothesis. The node factors are changed
as described in the proof of Lemma \ref{lem:7}. As a consequence,
when we sum over all possible renormalised self-energy clusters,
we find $\MM^{[\le n]}_{j}(x)=\MM^{[\le n]*}_{j}(x)$, which proves
the assertion.$\EP$

\begin{lemma}\label{lem:15}
Assume that the propagators $g^{[p]}_{j}(x)$ and their derivatives can
be uniformly bounded for all $0\le p \le n-1$ as in (\ref{eq:5.7})
and (\ref{eq:5.9}), for some constants $K_{1}$ and $K_{2}$.
Then one has $\MM^{[n]}_{1}(0)=-\MM^{[n]}_{2}(-2\la_{0})$
for all $n \ge 1$.
\end{lemma}

\prova Write $M^{[n]}_{1}(x)$ according to (\ref{eq:4.6}).
For any $T$ contributing to $M^{[n]}_{1}(x)$ we construct a renormalised
self-energy cluster $T'$ contributing to $M^{[n]}_{2}(x)$ as follows.
Call $\PP(T)$ the path of lines and nodes between the external
lines of $T$, and denote with $V(\PP(T))$ and $L(\PP(T))$ the set of
nodes and the set of lines, respectively, along $\PP(T)$. If $L(\PP(T))=
\emptyset$ the assertion trivially follows from (\ref{eq:3.23}).
Hence in the following assume $L(\PP(T)) \neq \emptyset$.

By definition of self-energy cluster all $\ell\in L(\PP(T))$
have momentum different from zero, while all lines connected to
a node $v\in V(\PP(T))$ have zero momentum (cf. Figure \ref{fig:3.8}).
Hence $\sum_{v\in V(\PP(T))}\nn_{v}=\vzero$. The nodes $v\in V(\PP(T))$
are totally ordered, so that we can number them $v_{0},v_{1},
\ldots,v_{N}$, if $N=|L(\PP(T))|$. The self-energy cluster
$T'$ is obtained through three steps: (i) first, we associate to each
node $v_{i}$, $i=0,\ldots,N$, the mode label and the node factor
of the node $v_{N-i}$ in $T$, -- in other words we revert the order of
the nodes, -- (ii) next, we write all node factors $f_{11,\nn_{v}}$
and $f_{22,\nn_{v}}$ as $f_{11,\nn_{v}}=-f_{22,\nn_{v}}$
and $f_{22,\nn_{v}}=-f_{11,\nn_{v}}$, -- by using that $\tr f =0$, --
(iii) finally we change consistently the component labels $j_{\ell}$
of the lines $\ell\in L(\PP(T))$, -- which means that each
label $j=1$ is changed into $j=2$ and vice versa.

If $\ell \in L(T)$ is the line connecting, say, $v_{k}$ to
$v_{k-1}$ for some $k=1,\ldots,N$, we still call $\ell$
the line in $L(T')$ which connects $v_{N-k+1}$ to $v_{N-k}$.
For each line $\ell \in L(T)$ we can write its momentum
as $\nn_{\ell}=\nn_{\ell}^{0}+\nn$, where $\nn_{\ell}^{0}$
is the sum of the mode labels of the nodes $v\in V(\PP(T))$
preceding $v$ in $T$ and $\nn$ is the momentum of the line entering $T$.
Then the corresponding line $\ell$ in $L(T')$ will have
momentum $-\nn_{\ell}^{0}+\nn$. Therefore each propagator
$g^{[n]}_{j}(\oo\cdot\nn_{\ell}+\oo\cdot\nn)$ in $T$ is changed
into $g^{[n]}_{3-j}(-\oo\cdot\nn_{\ell}+\oo\cdot\nn)$ in $T'$.

From the very definition of the propagators one sees immediately that,
by setting $x_{\ell}^{0}=\oo\cdot\nn_{\ell}^{0}$
and $x=\oo\cdot\nn$, one has
\begin{eqnarray}
& & g^{[n_{\ell}]}_{1}(x_{\ell}^{0}) =
g^{[n_{\ell}]}_{2} (x_{\ell}^{0}-2\la_{0}) =
- g^{[n_{\ell}]}_{2} ( - x_{\ell}^{0} - 2\la_{0}) , \nonumber \\
& & g^{[n_{\ell}]}_{2} (x_{\ell}^{0}) =
g^{[n_{\ell}]}_{1}(x_{\ell}^{0}+2\la_{0})  = 
- g^{[n_{\ell}]}_{1} (-x_{\ell}^{0} - 2\la_{0}) .
\label{eq:5.15} \end{eqnarray}

Now compute $\Val(T)$ for $x=0$ and $\Val(T')$ for $x=-2\la_{0}$.
Of course the node factors do not depend on the momenta, so that
\begin{equation}
\prod_{v\in V(\PP(T))} F_{v} =
(-1)^{|J(\PP(T))|} \prod_{v\in V(\PP(T'))} F_{v} ,
\label{eq:5.16} \end{equation}
where $J(\PP(T))$ is the set of nodes $v\in V(\PP(T))$ with
$j_{\ell_{v}}=j_{\ell_{v}'}$. It is immediate to realise that
$|J(\PP(T))|$ has the same parity of $|V(\PP(T))|$, -- 
see the proof of Lemma \ref{lem:7} for a similar argument.

By using (\ref{eq:5.15}) we obtain also
\begin{equation}
\prod_{\ell\in L(\PP(T))} g^{\RR}_{\ell} \Big|_{x=0} =
(-1)^{|L(\PP(T))|}
\prod_{\ell\in L(\PP(T'))} g^{\RR}_{\ell} \Big|_{x=-2\la_{0}} .
\label{eq:5.17} \end{equation}
Finally we have
\begin{equation}
\Big( \prod_{v\in P(T) \setminus V(\PP(T))} F_{v} \Big)
\Big( \prod_{\ell\in L(T) \setminus L(\PP(T))} g^{\RR}_{\ell} \Big)
= \Big( \prod_{v\in P(T') \setminus V(\PP(T'))} F_{v} \Big)
\Big( \prod_{\ell\in L(T') \setminus L(\PP(T'))} g^{\RR}_{\ell} \Big)
\label{eq:5.18} \end{equation}
for all $x\in\RRR$, so that, by using that
$(-1)^{|J(\PP(T))|}(-1)^{|L(\PP(T))|}=
(-1)^{|V(\PP(T))|+|L(\PP(T))|}=-1$, we find
$\VV_{T}(0) = - \VV_{T'}(-2\la_{0})$.
Then the assertion follows.$\EP$

\begin{lemma}\label{lem:16}
Assume that the propagators $g^{[p]}_{j}(x)$ and their derivatives can
be uniformly bounded for all $0\le p \le n-1$ as in (\ref{eq:5.7})
and (\ref{eq:5.9}), for some constants $K_{1}$ and $K_{2}$.
Then for $\eps$ small enough and $n \ge 1$ one has
$\MM^{[n]}_{1}(0)=\MM^{[n]}_{2}(-2\la_{0})=0$, and
\begin{eqnarray}
\left| \MM^{[n]}_{1}(x) \right| & \le &
B_{1} {\rm e}^{-\ka_{1} 2^{n}} |\eps|^{2}
\min\{ C_{1}^{-1} , |x| C_{1}^{-2}\} , \nonumber \\
\left| \MM^{[n]}_{2}(x) \right| & \le &
B_{1} {\rm e}^{-\ka_{1} 2^{n}} |\eps|^{2}
\min\{ C_{1}^{-1} , |x + 2 \la_{0}| C_{1}^{-2}\} ,
\label{eq:5.19} \end{eqnarray}
for suitable $n$-independent constants $B_{1}$ and $\ka_{1}$.
\end{lemma}

\prova By using the definitions in (\ref{eq:4.6}) and noting that all
sums are controlled, we see that the bound (\ref{eq:5.8})
implies the bound $| \MM^{[n]}_{j}(x) | \le
B_{1} {\rm e}^{-\ka_{1} 2^{n}} |\eps|^{2} C_{1}^{-1}$
for both $j=1$ and $j=2$.

The proof of the other bounds is more subtle.
Let us start with the case $j=1$.

Let $T$ be a renormalised self-energy cluster.
First consider the case $C_{1}\g_{n(M(T))} \le 4|\oo\cdot\nn|$,
where $\nn$ is the momentum associated with the entering line of $T$.
In that case one can extract from the last product in (\ref{eq:5.11})
a factor ${\rm e}^{-\ka_{0} M(T)/4} \le {\rm e}^{-\ka_{0} 2^{n(M(T))}/8}$.
Since $\oo$ is a Bryuno number then $a_{n}:=2^{-n}\log1/\al_{n}$
tends to zero as $n\to\infty$, hence for $\oo\cdot\nn$ small enough
one has ${\rm e}^{-\ka_{0} 2^{n(M(T))}/8} \le
(C_{0} \g_{n(M(T))})^{\ka_{0}/8a_{n(M(T))}} \le
C_{0} \g_{n(M(T))} \le 4 C_{0} C_{1}^{-1} |\oo\cdot\nn|$,
which implies the bound (\ref{eq:5.19}).

Then we consider the case $C_{1}n(M(T)) > 4|\oo\cdot\nn|$.
In that case for any line $\ell \in L(T)$ and for any $n<n_{\ell}$,
by the Diophantine conditions (\ref{eq:4.14}), one has
$|x_{\ell}^{0}+2\rho_{0}(x_{\ell}^{0})| >
C_{1}\g_{n(\nn_{\ell}^{0})}$, where $x_{\ell}=\oo\cdot\nn_{\ell}$ and
$x_{\ell}^{0}=\oo\cdot\nn_{\ell}^{0}$, with $\nn_{\ell}=
\nn_{\ell}^{0}+\sigma_{\ell}^{0}$, $\sigma_{\ell}=0,1$.
Since $|\nn_{\ell}^{0}| \le M(T)$, then
$C_{1}\g_{n(\nn_{\ell}^{0})} \ge C_{1}\g_{n(M(T))} > 4|\oo\cdot\nn|$,
which yields
\begin{equation}
2 \left| x_{\ell}^{0} + 2\rho_{0}(x_{\ell}^{0}) \right|
\ge \left| x_{\ell} + 2\rho_{0}(x_{\ell}) \right|
\ge \frac{1}{2}
\left| x_{\ell}^{0} + 2\rho_{0}(x_{\ell}^{0}) \right| .
\label{eq:5.20} \end{equation}
Such a property is important for the following reason. It can happen,
by the properties of the compact support functions, that a line $\ell$
is such that $g^{[n_{\ell}]}_{j_{\ell}}(x_{\ell}^{0}) \neq 0$,
whereas $g^{[n_{\ell}]}_{j_{\ell}}(x_{\ell}^{0}+x) =0$.
On the other hand in order to exploit the cancellations describe
below we have to consider also renormalised self-energy clusters
containing lines of this kind. Then (\ref{eq:5.20}) says that
in such cases, even if the bounds (\ref{eq:4.7}) are not satisfied,
one still has bounds of the same form with the only difference
that $\beta$ is replaced with $2\beta$ in the upper bound
and with $\beta/2$ in the lower bound. But this is enough
to apply both Lemma \ref{lem:11} and Lemma \ref{lem:12}.

For any  renormalised self-energy cluster we
consider the renormalised self-energy clusters
which belong to the same equivalence class.
Assume that $T$ is that of the second kind and that
$T'$ and $T''$ are those of the first type.
The corresponding self-energy values differ because of two facts:
(i) the value of $T'$ and $T''$ has an extra overall
factor $-1/2$, deriving from the product of the propagator $i$
times the node factor $i/2$, and (ii) for all lines along the path
between the external lines of $T$ the propagators depend also
on $\oo\cdot\nn$. The latter statement means that if $\ell$
is one of such lines then $g_{\ell}^{\RR}=g^{[n_{\ell}]}_{j_{\ell}}
(\oo\cdot\nn_{\ell}^{0}+\oo\cdot\nn)$ for $\ell\in L(T)$,
while $g_{\ell}^{\RR}=g^{[n_{\ell}]}_{j_{\ell}}
(\oo\cdot\nn_{\ell}^{0})$ for $\ell\in L(T')$ and $\ell\in L(T'')$. 
Finally, the two renormalised self-energy clusters
$T'$ and $T''$ have the same values.

Therefore $\VV_{T'}(\oo\cdot\nn)=\VV_{T'}(0)$ and
$\VV_{T''}(\oo\cdot\nn)=\VV_{T''}(0)=\VV_{T'}(0)$, hence
\begin{eqnarray}
& & \VV_{T'}(\oo\cdot\nn) + \VV_{T''}(\oo\cdot\nn) +
\VV_{T}(\oo\cdot\nn)
= \VV_{T'}(0) + \VV_{T''}(0) + \VV_{T}(0) \nonumber \\
& & \qquad \qquad + \left(
\VV_{T}(\oo\cdot\nn) - \VV_{T}(0) \right) =
\VV_{T}(\oo\cdot\nn) - \VV_{T}(0) ,
\label{eq:5.21} \end{eqnarray}
as $\VV_{T'}(0) = \VV_{T''}(0) = - \VV_{T}(0)/2$. By writing
\begin{equation}
\VV_{T}(\oo\cdot\nn) - \VV_{T}(0) =
\oo\cdot\nn \int_{0}^{1} {\rm d} s \, \frac{{\rm d}}{{\rm d}x}
\left. \VV_{T}(x) \right|_{x=s\oo\cdot\nn}
\label{eq:5.22} \end{equation}
and using (\ref{eq:5.10}) the bound (\ref{eq:5.19}) follows
once more.

The case $j=2$ follows easily from Lemma \ref{lem:15}.
Indeed for any renormalised self-energy cluster $T$ we can write
\begin{equation}
\VV_{T}(\oo\cdot\nn) = \VV_{T}(-2\la_{0}) +
\left( \VV_{T}(\oo\cdot\nn) - \VV_{T}(-2\la_{0}) \right) ,
\label{eq:5.23} \end{equation}
where
\begin{equation}
\VV_{T}(\oo\cdot\nn) - \VV_{T}(-2\la_{0}) =
\left( \oo\cdot\nn + 2\la_{0} \right)
\int_{0}^{1} {\rm d} s \, \frac{{\rm d}}{{\rm d}x}
\left. \VV_{T}(x) \right|_{x=-2\la_{0}+s(\oo\cdot\nn+2\la_{0})}
\label{eq:5.24} \end{equation}
can be bounded by using (\ref{eq:5.10}), while
\begin{equation}
\frac{i}{2} \sum_{k=1}^{\infty}
\sum_{T \in \SSS_{k,2,n-1}} \VV_{T}(-2\la_{0}) =
\MM^{[n]}_{2}(-2\la_{0}) = - \MM^{[n]}_{1}(0) = 0 ,
\label{eq:5.25} \end{equation}
so that the assertion is proved also in such a case.$\EP$

\begin{lemma}\label{lem:17}
Assume that the propagators $g^{[p]}_{j}(x)$ are differentiable,
and that, together with their derivatives, they
can be uniformly bounded for all $0\le p \le n-1$ as in (\ref{eq:5.7})
and (\ref{eq:5.9}), for suitable constants $K_{1}$ and $K_{2}$.
Then for $\eps$ small enough $\MM^{[\le n]}_{j}(x)$ is
differentiable in $x$, and one has
\begin{eqnarray}
& & \left| \MM^{[\le n]}_{j}(x') - \MM^{[\le n]}_{j}(x) -
\partial_{x} \MM^{[\le n]}_{j}(x) \left( x' - x \right) \right|
= o(\eps^{2} C_{1}^{-2}|x'-x|) , \nonumber \\
& & \left| \partial_{x} \MM^{[\le n]}_{j}(x) \right| \le
B_{2} |\eps|^{2} C_{1}^{-2} ,
\label{eq:5.26} \end{eqnarray}
for a suitable constant $B_{2}$.
\end{lemma}

\prova By writing $\MM^{[\le n]}_{j}(x)$ according to (\ref{eq:4.6}),
one finds immediately that the function is differentiable
if the propagators are differentiable, and that
the derivative satisfies the bound in (\ref{eq:5.19}).
The factor $\eps^{2}$ is due to the fact that a self-energy cluster
$T$ depending explicitly on $x$ has at least $k_{T}=2$.$\EP$

\begin{lemma}\label{lem:18}
Assume that the propagators $g^{[p]}_{j}(x)$ and their derivatives
can be uniformly bounded for all $0\le p \le n-1$ as in (\ref{eq:5.7})
and (\ref{eq:5.9}), for some constants $K_{1}$ and $K_{2}$.
Then for $\eps$ small enough one has
\begin{equation}
\left| x + 2 \MM^{[\le n]}_{j}(x) \right| \ge
\frac{1}{2} \Delta_{0}(x) 
\label{eq:5.27} \end{equation}
as far as $g^{[n]}_{j}(x) \neq 0$.
\end{lemma}

\prova By Lemma \ref{lem:16} one has $\MM^{[\le n]}_{1}(0)=0$
and $\MM^{[\le n]}_{2}(-2\la_{0})=\la_{0}$. Set $j(x)=1$
when $\rho_{0}(x)=0$ and $j(x)=2$ when $\rho_{0}(x)=\la_{0}$,
so that one can write
\begin{eqnarray}
x + 2\MM^{[\le n]}_{j(x)}(x) & = &
x + 2\MM^{[\le n]}_{j(x)}(-2\rho_{0}(x)) +
\left( 2\MM^{[\le n]}_{j(x)}(x) -
2\MM^{[\le n]}_{j(x)}(-2\rho_{0}(x)) \right)
\nonumber \\
& = & x + 2 \rho_{0}(x) + 2 \left( \MM^{[\le n]}_{j(x)}(x) -
\MM^{[\le n]}_{j(x)}(-2\rho_{0}) \right) ,
\label{eq:5.28} \end{eqnarray}
where $| \MM^{[\le n]}_{j(x)}(x) - \MM^{[\le n]}_{j(x)}(-2\rho_{0})|
\le \hbox{const.} |\eps|^{2} C_{1}^{-2} |x + 2\rho_{0}(x)|$,
by Lemma \ref{lem:17}. Then by (\ref{eq:4.7}) one has
\begin{equation}
\left| x + 2 \MM^{[\le n]}_{j(x)}(x) \right| \ge
\left(1 - \hbox{const.} |\eps|^{2} C_{1}^{-2} \right)
\left| x+2\rho_{0}(x) \right| \ge
\frac{1 - \hbox{const.} |\eps|^{2} C_{1}^{-2}}{\sqrt{2}}
\Delta_{0}(x) .
\label{eq:5.29} \end{equation}
Since $| x + 2 \MM^{[\le n]}_{3-j(x)}(x) | \ge
(1-\hbox{const.}|\eps|^{2}C_{1}^{-2}) | x + 2 \MM^{[\le n]}_{j(x)}(x)|$,
the bound follows.$\EP$

\begin{lemma}\label{lem:19}
The propagators $g^{[n]}_{j}(x)$ satisfy the bounds (\ref{eq:5.7})
and (\ref{eq:5.9}) for all $n\ge 0$.
\end{lemma}

\prova The proof can be performed by induction.
For $n=1$ the bounds (\ref{eq:5.7}) and (\ref{eq:5.9})
are trivially satisfied, as $\MM^{[0]}_{1}(x)=0$ and
$\MM^{[0]}_{2}(x)=\la_{0}$, because of the
Diophantine conditions (\ref{eq:4.14}).

The difference for $n>1$ is that now the propagators depend
also on the functions $\MM^{[p]}_{j}(x)$, $p<n$,
appearing in the denominators and the compact support functions.
Then assume (\ref{eq:5.7}) and (\ref{eq:5.9})
for all $p<n$. Then one has
$|g^{[n]}_{j}(x)| \le \hbox{const.} \psi_{n}(\Delta_{0}(x))/
\Delta_{0}(x) \le \hbox{const.} C_{1}^{-1} \g_{n}^{-1}$,
by Lemma \ref{lem:18}. Moreover
\begin{eqnarray}
\partial_{x} g^{[n]}_{j}(x) & = & -i
\sum_{p=0}^{n-1} \chi_{0}(\Delta_{0}(x)) \ldots
\partial \chi_{p}(\Delta_{0}(x)) \ldots \chi_{n-1}(\Delta_{0}(x))
\psi_{n}(\Delta_{0}(x)) \frac{\partial_{x}\Delta_{0}(x)}
{x+2 \MM^{[\le n]}_{j}(x)} \nonumber \\
& - &  i \chi_{0}(\Delta_{0}(x)) \ldots \chi_{n-1}(\Delta_{0}(x))
\partial \psi_{n}(\Delta_{0}(x)) 
\frac{\partial_{x}\Delta_{0}(x)}
{x+2 \MM^{[\le n]}_{j}(x)} \nonumber \\
& + & i \chi_{0}(\Delta_{0}(x)) \ldots \chi_{n-1}(\Delta_{0}(x))
\psi_{n}(\Delta^{[n]}(x)) \frac{1+2\partial_{x}\MM^{[\le n]}_{j}(x)}
{(x+2 \MM^{[\le n]}_{j}(x))^{2}} ,
\label{eq:5.30} \end{eqnarray}
where $\partial$ denotes derivative with respect to the argument.

One checks immediately that for all $p \ge 0$
\begin{equation}
\partial \chi_{p}(x) \le \hbox{const.} C_{1}^{-1} \g_{p}^{-1} ,
\qquad \partial \psi_{p}(x) \le \hbox{const.} C_{1}^{-1} \g_{p}^{-1} ,
\qquad \partial_{x} \Delta_{0}(x) \le \hbox{const.}  ,
\label{eq:5.31} \end{equation}
so that the derivative $\partial_{x}\MM^{[\le n]}_{j}(x)$ can be bounded
through (\ref{eq:5.26}), because of the inductive hypothesis.

Hence, by using once more (\ref{eq:4.9}) and Lemma \ref{lem:18}
to bound the denominators, we obtain from (\ref{eq:5.30})
\begin{equation}
\left| \partial_{x} g^{[n']}_{j}(x) \right| \le
\hbox{const.} C_{1}^{-2} \left( \sum_{p=0}^{n-1} \g_{p}^{-1}
\g_{n}^{-1} + \g_{n}^{-1} \g_{n}^{-1} + \g_{n}^{-2} \right) \le
\hbox{const.} C_{1}^{-2} \g_{n}^{-3} ,
\label{eq:5.32} \end{equation}
which proves the assertion.$\EP$

\begin{lemma}\label{lem:20}
Let $\la\in\Lambda_{0}^{*}$. There exists $\eps_{0}>0$ such that
for $|\eps|<\eps_{0}$ the coefficients $u^{[k]}_{j,\nn}$, $j=1,2$, and
$\mu^{[k]}$ are bounded by
\begin{equation}
\left| u^{[k]}_{j,\nn} \right| \le B \, {\rm e}^{-\ka|\nn|}|\eps|^{k} ,
\qquad \left| \mu^{[k]} \right| \le B \, |\eps|^{k} ,
\label{eq:5.33} \end{equation}
for suitable $k$-independent constants $B$ and $\ka$.
One can take $\eps_{0}=O(C_{1} \g_{m_{0}})$, with $m_{0}$
depending on $\ka_{0}$.
\end{lemma}

\prova For any tree $\theta\in\Theta^{\RR}_{k,j,\nn}$
the value $\Val(\theta)$ can be bounded by using the bounds
(\ref{eq:5.11}) for the factors $F_{v}$ and the bounds
(\ref{eq:5.7}), proved in Lemma \ref{lem:19}, for the propagators.
Summation over the Fourier labels can be performed by
using an exponential decay factor ${\rm e}^{-\ka_{0}M(T)/4}$
which can be extracted from (\ref{eq:5.11}). Summation
over the other labels and over the number of unlabelled trees
can be easily bounded as a constant to the power $k$.$\EP$

\begin{lemma}\label{lem:21}
The function $\overline u(t)$ solves (\ref{eq:3.3})
for all $\nn\neq\vzero$, provided $\mu=\overline \mu$.
\end{lemma}

\prova We write
\begin{eqnarray}
\overline u_{j}(t) & = & \overline u_{j,\vzero} +
\sum_{\nn \in\ZZZ^{d}} e^{i\oo\cdot\nn} \overline u_{j,\nn} ,
\qquad \overline u_{j,\nn} = \sum_{n=0}^{\infty}
\overline u_{j,\nn,n} , \nonumber \\
\overline u_{j,\nn,n} & = &
\sum_{k=1}^{\infty} \eps^{k} \sum_{\theta\in\Theta^{\RR}_{k,j,\nn,n}}
\Val(\theta) ,
\label{eq:5.34} \end{eqnarray}
where $\Theta^{\RR}_{k,j,\nn,n}$ is the set of trees
in $\Theta^{\RR}_{k,j,\nn}$ with root line on scale $n$.

An important property of the compact support functions is that
\begin{equation}
1 = \sum_{n=0}^{\infty} \Psi_{n}(x) , \qquad
\Psi_{n}(x) := \chi_{0}(\Delta_{0}(x)) \ldots
\chi_{n-1} (\Delta_{0}(x)) \psi_{n} (\Delta_{0}(x)) ,
\label{eq:5.35} \end{equation}
where the summand for $n=0$ is meant as $\psi_{0}(\Delta_{0}(x))$.
More generally one has for all $s\ge 1$
\begin{equation}
1 = \sum_{n=p}^{\infty} \Psi_{p,n}(x) , \qquad
\Psi_{p,n}(x) := \chi_{p}(\Delta_{0}(x)) \ldots
\chi_{n-1} (\Delta_{0}(x)) \psi_{n} (\Delta_{0}(x)) ,
\label{eq:5.36} \end{equation}
where again the summand for $n=p$ is meant as $\psi_{p}(\Delta_{0}(x))$.

We can rewrite the equation (\ref{eq:3.3}) as
\begin{equation}
u_{j,\nn} = g_{j} (x) \Phi_{j,\nn}(u) ,
\qquad \Phi_{j} = \eps f_{j1} + i \mu + \eps f_{j1} u_{1} +
\eps f_{j2} u_{2} + (-1)^{j+1} i \mu u_{j} , 
\label{eq:5.37} \end{equation}
where $x=\oo\cdot\nn$, $g_{j}(x)=-i(x+2\MM^{[0]}_{j}(x))^{-1}$,
with $\MM^{[0]}_{1}(x)=0$ and $\MM^{[0]}_{2}(x)=\la_{0}$.

By using (\ref{eq:5.35}) we can write
\begin{eqnarray}
g_{j} (x) \Phi_{j,\nn}(\overline u) & = &
g_{j}(x) \sum_{n=0}^{\infty} \Psi_{n}(x) \Phi_{j,\nn}(\overline u)
\nonumber \\
& = & g_{j}(x) \sum_{n=0}^{\infty} \Psi_{n}(x)
\left( g^{[n]}_{j}(x) \right)^{-1}
\left( g^{[n]}_{j}(x) \Phi_{j,\nn}(\overline u) \right) ,
\label{eq:5.38} \end{eqnarray}
where $\Psi_{n}(x) ( g^{[n]}_{j}(x) )^{-1} =
i ( x + 2\MM^{[\le n]}_{j}(x) )$, and
\begin{equation}
g^{[n]}_{j}(x) \Phi_{j,\nn}(\overline u) =
\sum_{k=1}^{\infty} \eps^{k}
\sum_{\theta\in\overline\Theta^{\RR}_{k,j,\nn,n}}
\Val(\theta) ,
\label{eq:5.39} \end{equation}
where $\overline\Theta^{\RR}_{k,j,\nn,n}$
differs from $\Theta^{\RR}_{k,j,\nn,n}$ as it contains also
trees which can have one renormalised self-energy cluster $T$
with exiting line given by the root line of $\theta$.
In such a case if $p$ is the line of the entering line of $T$,
then $p\ge 0$ and the scale $n_{T}$ of $T$ is such that
$n_{T}+1 \le \min\{n,p\}$, by definition of cluster.

Then we have
\begin{eqnarray}
& & \sum_{n=0}^{\infty} \Psi_{n}(x) \left( g^{[n]}_{j}(x) \right)^{-1}
\left( g^{[n]}_{j}(x) \Phi_{j,\nn}(u) \right) \nonumber \\
& & \qquad = i 
\sum_{n=0}^{\infty} \left( x + 2\MM^{[\le n]}_{j}(x) \right) 
\sum_{k=1}^{\infty} \eps^{k}
\sum_{\theta\in\Theta^{\RR}_{k,j,\nn,n}} \Val(\theta) \nonumber \\
& & \qquad - 2 i
\sum_{n=1}^{\infty}
\Psi_{n}(x) \sum_{p=n}^{\infty} \sum_{s=1}^{n} M^{[s]}_{j}(x)
\sum_{k=1}^{\infty} \eps^{k}
\sum_{\theta\in\Theta^{\RR}_{k,j,\nn,p}} \Val(\theta) \nonumber \\
& & \qquad - 2 i
\sum_{n=2}^{\infty}
\Psi_{n}(x) \sum_{p=1}^{n-1} \sum_{s=1}^{p} M^{[s]}_{j}(x)
\sum_{k=1}^{\infty} \eps^{k}
\sum_{\theta\in\Theta^{\RR}_{k,j,\nn,p}} \Val(\theta) ,
\label{eq:5.40} \end{eqnarray}
and we can use the definitions (\ref{eq:5.34}) to write
\begin{equation}
\sum_{k=1}^{\infty} \eps^{k}
\sum_{\theta\in\Theta^{\RR}_{k,j,\nn,n}} \Val(\theta)
= \overline u_{j,\nn,n} , \qquad
\sum_{k=1}^{\infty} \eps^{k}
\sum_{\theta\in\Theta^{\RR}_{k,j,\nn,p}} \Val(\theta)
= \overline u_{j,\nn,p} ,
\label{eq:5.41} \end{equation}
in the second line and, respectively, in the third and fourth lines.

Then the sum of the third and fourth lines in (\ref{eq:5.40}) gives
\begin{eqnarray}
& & -2i \left( \sum_{n=1}^{\infty} \Psi_{n}(x)
\sum_{p=n}^{\infty} \sum_{s=1}^{n} M^{[s]}_{j}(x)
\, \overline u_{j,\nn,p} +
\sum_{n=2}^{\infty} \Psi_{n}(x)
\sum_{p=1}^{n-1} \sum_{s=1}^{p} M^{[s]}_{j}(x)
\, \overline u_{j,\nn,p} \right) \nonumber \\
& & \qquad \qquad = -2i
\sum_{p=1}^{\infty} \overline u_{j,\nn,p} \left(
\sum_{s=1}^{p} \sum_{n=s}^{p} M^{[s]}_{j}(x) \Psi_{n}(x) +
\sum_{s=1}^{p} \sum_{n=p+1}^{\infty} M^{[s]}_{j}(x) \Psi_{n}(x)
\right) \nonumber \\
& & \qquad \qquad = -2i
\sum_{p=1}^{\infty} \overline u_{j,\nn,p} \sum_{s=1}^{p}
M^{[s]}_{j}(x) \sum_{n=s}^{\infty} \Psi_{n}(x) .
\label{eq:5.42} \end{eqnarray}
If we define
\begin{equation}
\Xi_{n}(x) := \chi_{0} (\Delta_{0}(x)) \ldots
\chi_{n-1} (\Delta_{0}(x)) \chi_{n} (\Delta_{0}(x)) ,
\label{eq:5.43} \end{equation}
then in (\ref{eq:5.42}) we can write
\begin{equation}
\sum_{n=s}^{\infty} \Psi_{n}(x) =
\Xi_{s-1}(x) \sum_{q=s}^{\infty} \Psi_{s,q} = \Xi_{s-1}(x) ,
\label{eq:5.44} \end{equation}
where the property (\ref{eq:5.36}) has been used.
Hence in (\ref{eq:5.42}) we have
\begin{equation}
\sum_{s=1}^{p} M^{[s]}_{j}(x) \sum_{n=s}^{\infty} \Psi_{n}(x) =
\sum_{s=1}^{p} M^{[s]}_{j}(x) \Xi_{s-1}(x) =
\MM^{[\le p]}_{j}(x) - \MM^{[0]}_{j}(x) ,
\label{eq:5.45} \end{equation}
where the factor $\MM^{[0]}_{j}(x)$ has been subtracted
as the sum over $s$ starts from $s=1$ and not from $s=0$.

If we insert (\ref{eq:5.40}) into (\ref{eq:5.38}),
by taking into account (\ref{eq:5.42}) and (\ref{eq:5.45}),
we obtain
\begin{eqnarray}
g_{j}(x) \Phi_{j,\nn}(\overline u) & = & g_{j}(x) \left( i
\sum_{n=0}^{\infty} \left( x + 2\MM^{[\le n]}_{j}(x) \right) -
2 i \sum_{n=1}^{\infty} \left( \MM^{[\le n]}_{j}(x) -
\MM^{[0]}_{j}(x) \right) \right) \overline u_{j,\nn,n}
\nonumber \\
& = & g_{j}(x) \sum_{n=0}^{\infty} i \left( x + 2\MM^{[0]}_{j}(x) \right)
\overline u_{j,\nn,n} = \sum_{n=0}^{\infty}
\overline u_{j,\nn,n} = \overline u_{j,\nn} ,
\label{eq:5.46} \end{eqnarray}
so that (\ref{eq:5.37}) is satisfied for $u=\overline u$.$\EP$

\begin{lemma}\label{lem:22}
The function $\overline u(t)$ solves the system of differential
equations (\ref{eq:2.19}) for all $t\in\RRR$,
provided $\mu=\overline \mu$. Moreover
the function $H$ defined in (\ref{eq:2.20}) satisfies
$H(\overline u(t)) =0$ for all $t\in\RRR$.
Both $\overline u(t)$ and $\overline \mu$ are analytic in $\eps$.
\end{lemma}

\prova Because of Lemma \ref{lem:21}, to show that $\overline u(t)$
is a solution it is enough to prove that $\overline u$
solves (\ref{eq:3.4}), that is $0 = \Phi_{j,\vzero}(\overline u)$,
with $\Phi_{j}(u)$ defined in (\ref{eq:5.33}).
But this is obvious by construction.

The claim on $H(\overline u(t))$ follows if the solution
is in $\gotM$, so that (\ref{eq:2.16}) is satisfied.
But again this follows from the construction of the solution.

Finally the statement about analyticity easily follows
from the construction of the renormalised series.
The series defining $\overline u(t)$ in (\ref{eq:4.12})
and $\overline \mu$ in (\ref{eq:4.11}) can be viewed
as power series in $\eps$ with coefficients depending on $\eps$.
The coefficients depend on $\eps$ through the propagators,
and in fact are analytic in $\eps$ (for $\eps$ small enough).
Hence the series themselves define functions which are analytic
in $\eps$.$\EP$

A result analogous to Lemma \ref{lem:22}, in particular analyticity
of the conjugation and of the counterterm, was proved
by using renormalisation group techniques in \cite{L,KL}.
In the case of the Schr\"odinger equation it was also obtained
in \cite{BGGM}, with techniques similar to those used in this paper;
cf. also \cite{G0,E1,GM,BG} for related issues. See also \cite{GBG},
Chapter 9, for resummed series defining analytic functions,
in the context of maximal KAM tori.

\zerarcounters
\section{Reducibility on a large measure set}
\label{sec:6}

So far we have proved that, as far as $\la_{0}\in\Lambda_{0}^{*}$,
the function $\overline u(t)$ solves (\ref{eq:2.19}).
We still have to prove that the relative measure of $\Lambda_{0}^{*}$
with respect to $\Lambda_{0}$ is large, and we have to see
what this means for the parameter $\la\in[a,b]$. We shall find that
the subset $\Lambda^{*}\subset[a,b]$ of values $\la$ for which the
construction envisaged in the previous sections works is of large
measure; this will complete the proof of Theorem \ref{thm:1}.

As a consequence of Lemma \ref{lem:19},
we have that (\ref{eq:5.19}), (\ref{eq:5.26}),
and (\ref{eq:5.27}) hold for all $n \ge 0$.

For each $\nn\in\ZZZ^{d}_{*}$ we have to exclude all values $\la_{0}
\in\Lambda_{0}^{*}$ such that $|\oo\cdot\nn + 2\la_{0}|\le
C_{1}\g_{n(\nn)}$. If we consider $\la_{0}$ as a function of an
auxiliary parameter $t\in[-1,1]$, we can write
\begin{equation}
\oo\cdot\nn + 2\la_{0}(t) = t C_{1} \g_{n(\nn)} ,
\qquad t\in[-1,1] ,
\label{eq:6.1} \end{equation}
so that
\begin{equation}
\frac{{\rm d}\la_{0}}{{\rm d}t} = \frac{C_{1}}{2} \g_{n(\nn)} .
\label{eq:6.2} \end{equation}
Then for each $\nn\in\ZZZ^{d}_{*}$ we have to exclude
all values of $t$ in $[-1,1]$.

\begin{lemma}\label{lem:23}
There exists $\eps_{0}>0$ and $\sigma>0$ such that for
all $|\eps|<\eps_{0}$ the Lebesgue measure
of the set $\Lambda_{0}\setminus\Lambda_{0}^{*}$
is bounded proportionally to $|\eps|^{\sigma}$.
\end{lemma}

\prova The set $\Lambda^{*}_{0}$ is obtained by imposing
the Diophantine conditions (\ref{eq:4.14}). Then we can bound
\begin{equation}
\hbox{meas}(\Lambda_{0}^{*}) =
\int_{\Lambda_{0}^{*}} {\rm d}\lambda_{0} =
\sum_{\nn\in\ZZZ^{d}_{*}}
\int_{-1}^{1} {\rm d} t \, \frac{{\rm d}\lambda_{0}}{{\rm d}t}
= \sum_{\nn\in\ZZZ^{d}_{*}} C_{1} \g_{n(\nn)} ,
\label{eq:6.3} \end{equation}
where we can write
\begin{equation}
\sum_{\nn\in\ZZZ^{d}_{*}} \g_{n(\nn)} =
\sum_{n=0}^{\infty} \sum_{2^{n-1}<|\nn| \le 2^{n}}
\g_{n} \le \hbox{const.} \sum_{n=0}^{\infty} 2^{n(d-1)} \g_{n}
\le \hbox{const.}
\label{eq:6.4} \end{equation}
By inserting (\ref{eq:6.4}) into (\ref{eq:6.3}) we obtain
$\hbox{meas}(\Lambda_{0}^{*}) \le \hbox{const.} C_{1}$,
hence $\Lambda_{0}^{*}$ is a set of measure proportional to $C_{1}$.

By Lemma \ref{lem:20} we can take $C_{1}=|\eps|^{\sigma}$,
with $0<\sigma<1$. Hence the measure of the discarded set
can be bounded proportionally to $|\eps|^{\sigma}$.$\EP$

In the following write $\overline u(t,\la_{0})$, $\mu(\la_{0})$ and
$\MM^{[\le n]}_{j}(x,\la_{0})$ to make explicit the dependence of
$\overline u(t)$, $\mu$ and $\MM^{[\le n]}_{j}(x)$ on $\la_{0}$.

\begin{lemma}\label{lem:24}
Assume $\la_{0},\la_{0}'\in\Lambda_{0}^{*}$. One has
\begin{eqnarray}
& & \left| \MM^{[\le n]}_{j}(x,\la_{0}') -
\MM^{[\le n]}_{j}(x,\la_{0}) - \partial_{\la_{0}}
\MM^{[\le n]}_{j}(x,\la_{0}) \right| =
o(\eps^{2}C_{1}^{-2} |\la_{0}'-\la_{0}|) ,
\nonumber \\
& & \left| \partial_{\la_{0}} \left(
\MM^{[\le n]}_{j}(x,\la_{0}) - \MM^{[0]}_{j}(x,\la_{0}) \right)
\right| \le A_{1} |\eps|^{2} C_{1}^{-2} ,
\label{eq:6.5} \end{eqnarray}
for a suitable constant $A_{1}$. In particular $\MM^{[\le n]}_{j}
(x,\la_{0})$ can be extended in all $\overline \Lambda_{0}$
to a differentiable function (Whitney extension),
whose derivative satisfies the bound in (\ref{eq:6.5}).
\end{lemma}

\prova The proof is by induction. $\MM^{[\le n]}_{j}(x,\la_{0})$
can be written according to (\ref{eq:4.6}). Hence it depends recursively
on $\MM^{[\le n']}_{j'}(x)$, $n'<n$, through the propagators,
and one can express $\partial_{\la_{0}} \MM^{[\le n]}_{j}(x,\la_{0})$
as sum of derivatives of self-energy values,
\begin{equation}
\partial_{\la_{0}} \VV_{T}(x) = \eps^{k_{T}}
\sum_{\ell\in \PP(T)} \partial_{\la_{0}} g_{\ell}^{\RR}
\Big( \prod_{\ell'\in L(T) \setminus \ell} g_{\ell'}^{\RR} \Big)
\Big( \prod_{v \in P(T)} F_{v} \Big) .
\label{eq:6.6} \end{equation}

For $n=0$ the assertion is trivially satisfied, as
$g^{[0]}_{j}(x)=-i(x+2\MM^{[0]}_{j}(x,\la_{0}))^{-1}$, with
$\MM^{[0]}_{1}(x,\la_{0}))=0$ and $\MM^{[0]}_{2}(x,\la_{0}))=\la_{0}$.
Then for all $\la_{0}\in \Lambda_{0}^{*}$ one has
$\partial_{\la_{0}} g^{[0]}_{1}(x,\la_{0})=0$
and $\partial_{\la_{0}} g^{[0]}_{2}(x,\la_{0})=-2i
(x+2\la_{0})^{-2}$, and one can consider the Whitney extension of
$\MM^{[0]}_{j}(x,\la_{0}))$ in all $\overline\Lambda_{0}$.

For $n \ge 1$ assume that all $\MM^{[\le n']}_{j'}(x,\la_{0})$, $n'<n$,
can be extended to differentiable functions in $\overline\Lambda_{0}$
and satisfy the bounds in (\ref{eq:6.5}).
Then the derivative $\partial_{\la_{0}} g_{\ell}^{\RR}$,
in (\ref{eq:6.6}), can be bounded because of the inductive hypothesis.
Simply, one reasons as in the proof of Lemma \ref{lem:19},
and (\ref{eq:6.5}) follows.$\EP$

The Whitney differentiability of $\MM^{[\le n]}_{j}(x,\la_{0})$
implies also that of $\overline u(t,\la_{0})$ and $\mu(\la_{0})$.
Hence the following result follows immediately from Lemma \ref{lem:24}.

\begin{lemma}\label{lem:25}
The renormalised series for $\overline u(t)$ and $\mu$
converge to functions which are differentiable
in the sense of Whitney in $\Lambda_{0}^{*}$.
\end{lemma}

Now, we can conclude the proof of theorem \ref{thm:1}.

\begin{lemma}\label{lem:26}
Call $\Lambda^{*}$ the subset of $[a,b]$ for which
the system (\ref{eq:2.19}) is reducible. There exists $\eps_{0}>0$
such that for all $|\eps|<\eps_{0}$ the Lebesgue measure
of the set $[a,b]\setminus\Lambda^{*}$ is bounded
proportionally to $|\eps|^{\sigma}$, with
$\sigma$ as in Lemma \ref{lem:23}.
\end{lemma}

\prova Write $\la_{0}+\mu=\la$. We want to fix the set $\Lambda_{0}$
so that for $\la_{0}\in\Lambda_{0}^{*} \subset \Lambda_{0}$
one has $\la \in \Lambda:=[a,b]$. Write $\Lambda_{0}=[a_{0},b_{0}]$,
with $a_{0} = a - \eps f_{11,\vzero} + A \eps^{2} C_{1}^{-1}$ and
$b_{0} = b - \eps f_{11,\vzero} - A \eps^{2} C_{1}^{-1}$,
where $A$ is a constant such that for all $\la_{0}\in\Lambda_{0}^{*}$
and all $|\eps|<\eps_{0}$ one has $|\mu-\eps f_{11,\vzero}| < A \eps^{2}
C_{1}^{-1}$ (this is possible by Lemma \ref{lem:20}). Then
$\hbox{meas}(\Lambda_{0})=\hbox{meas}(\Lambda)-2A\eps^{2}C_{1}^{-1}$,
whereas $\hbox{meas}(\Lambda_{0}^{*})=\hbox{meas}(\Lambda_{0})-
O(|\eps|^{\sigma})$
by Lemma \ref{lem:23}. Call $\Lambda^{*}$ the subset of values
$\la \in \Lambda$ such that $\la=\la_{0} + \mu$,
for $\la\in\Lambda_{0}^{*}$ and $\mu=\mu(\la_{0})$.
By construction $\Lambda^{*}\subset\Lambda$.

By Lemma \ref{lem:25} the function $\la_{0} \to \mu(\la_{0})$ is
differentiable in the sense of Whitney, so that
\begin{equation}
\frac{{\rm d}\la}{{\rm d}\la_{0}} =
\frac{{\rm d}}{{\rm d}\la_{0}} \left( \la_{0} + \mu \right)
= 1 + \frac{{\rm d}\mu}{{\rm d}\la_{0}} , \qquad
\left| \frac{{\rm d}\mu}{{\rm d}\la_{0}} \right|
\le \hbox{const.} |\eps|^{2} C_{1}^{-1} ,
\label{eq:6.7} \end{equation}
where we explicitly used that the first contribution
to $\mu$ depending on $\la_{0}$ has size $O(\eps^{2} C_{1}^{-1})$.

Therefore
\begin{equation}
\hbox{meas}(\Lambda\setminus\Lambda^{*})
= \int_{\Lambda\setminus\Lambda^{*}} {\rm d}\la \le
-2A \eps^{2} + \int_{\Lambda_{0}\setminus\Lambda^{*}_{0}} {\rm d}\la_{0}
\left| \frac{{\rm d}\la}{{\rm d}\la_{0}} \right|
\le \hbox{const.} |\eps|^{\sigma} ,
\label{eq:6.8} \end{equation}
because $\hbox{meas}(\Lambda_{0}\setminus\Lambda^{*}_{0})=
O(|\eps|^{\sigma})$ by Lemma \ref{lem:23},
and the assertion is proved.$\EP$

So far we assumed $0\notin[a,b]$. If $0\in[a,b]$ we can discard
a subset $\Lambda_{1}\subset[a,b]$ around the origin, of measure
$O(|\eps|^{\sigma})$, such that for all $\la\in[a,b] \setminus
\Lambda_{1}$ one has $|\la|>\hbox{const.}|\eps|^{\sigma}$. Then
$|\la_{0}|$ is bounded below proportionally to $|\eps|^{\sigma}$,
because $|\la-\la_{0}| =|\mu|=O(|\eps|)$ and $\sigma<1$.
Though, this does not modify the bounds of the previous sections.
Indeed the only difference is that also the propagators
with vanishing momentum (that is on scale $-1$) are
bounded proportionally to $|\eps|^{-\sigma}$ -- like
those with non-zero momentum $\nn_{\ell}$, which are bounded
proportionally to $|\eps|^{-\sigma}\g_{n(\nn_{\ell})}^{-1}$ --
and the bounds were obtained
by using that one has at worst a factor $|\eps|^{-\sigma}$ per line.
Then one can restrict the analysis to $[a,b]\setminus\Lambda_{1}$,
and the same conclusions hold.

\vskip1.truecm

\noindent \textbf{Acknowledgments.} I'm indebted to Giovanni
Gallavotti for useful discussions.



\begin{thebibliography}{99}

\bibitem{AK}{
A. Avila, R. Krikorian,
\textit{Reducibility or non-uniform hyperbolicity for quasi-periodic
Schr\"odinger cocycles},
Preprint, to appear on Ann. Math.}
%
\bibitem{BG}{
M. Bartuccelli, G. Gentile,
\textit{Lindstedt series for perturbations of isochronous systems},
Rev. Math. Phys.
\textbf{14} (2002), no. 2, 121--171. }
%
\bibitem{BGGM}{
F. Bonetto, G. Gallavotti, G. Gentile, V. Mastropietro,
\textit{Quasi-linear flows on tori: regularity of their linearization},
Comm. Math. Phys.
\textbf{192} (1998),  no. 3, 707--736.}
%
\bibitem{Ch}{
Ch.-Q. Cheng,
\textit{Birkhoff-Kolmogorov-Arnold-Moser tori in convex Hamiltonian systems},
Comm. Math. Phys
\textbf{177} (1996), no. 3, 529--559.}
%
\bibitem{DS}{
E.I. Dinaburg, Ja.G. Sina\u\i,
\textit{The one-dimensional Schr\"odinger equation
with quasiperiodic potential},
Funkcional. Anal. i Prilo\v zen.
\textbf{9} (1975), no. 4, 8--21;
English translation:
Functional Anal. Appl.
\textbf{9} (1975), no. 4, 279--289 (1976).}
%
\bibitem{E1}{
L.H. Eliasson,
\textit{Hamiltonian systems with linear normal form near
an invariant torus},
Nonlinear dynamics (Bologna, 1988), 11--29, World Sci. Publishing,
Teaneck, NJ, 1989.}
%
\bibitem{E2}{
L.H. Eliasson,
\textit{Floquet solutions for the $1$-dimensional quasi-periodic
Schr\"odinger equation},
Comm. Math. Phys.
\textbf{146} (1992), no. 3, 447--482.}
%
\bibitem{E3}{
L.H. Eliasson,
\textit{Reducibility and point spectrum for
linear quasi-periodic skew-products},
Proceedings of the International Congress of Mathematicians,
Vol. II (Berlin, 1998), Doc. Math. 1998, Extra Vol. II, 779--787.}
%
\bibitem{E4}{
L.H. Eliasson,
\textit{On the discrete one-dimensional quasi-periodic Schr\"odinger
equation and other smooth quasi-periodic skew products},
Hamiltonian systems with three or more degrees of freedom (S'Agaró, 1995),
55--61, NATO Adv. Sci. Inst. Ser. C Math. Phys. Sci., 533, 
Kluwer Acad. Publ., Dordrecht, 1999.}
%
\bibitem{G0}{
G. Gallavotti,
\textit{A criterion of integrability for perturbed nonresonant
harmonic oscillators. "Wick ordering" of the perturbations in classical
mechanics and invariance of the frequency spectrum},
Comm. Math. Phys.
\textbf{87} (1982/83), no. 3, 365--383.}
%
\bibitem{GBG}{
G. Gallavotti, F. Bonetto, G. Gentile,
\textit{Aspects of ergodic, qualitative and statistical theory of motion},
Texts and Monographs in Physics, Springer, Berlin, 2004.}
%
\bibitem{GGG}{
G. Gallavotti, G. Gentile, A. Giuliani,
\textit{Fractional Lindstedt series},
Preprint, to appear on J. Math. Phys.}
%
\bibitem{GG1}{
G. Gallavotti, G. Gentile,
\textit{Hyperbolic low-dimensional invariant tori
and summations of divergent series},
Comm. Math. Phys.
\textbf{227} (2002), no. 3, 421--460.}
%
\bibitem{G1}{
G. Gentile,
\textit{Quasi-periodic solutions for two-level systems},
Comm. Math. Phys.
\textbf{242} (2003), no. 1, 221--250.}
%
\bibitem{G2}{
G. Gentile,
\textit{Degenerate lower-dimensional tori under the Bryuno condition},
Preprint.}
%
\bibitem{GBC}{
G. Gentile, D.A. Cortez, J.C.A. Barata,
\textit{Stability for quasi-periodically perturbed Hill's equations},
Comm. Math. Phys.
\textbf{260} (2005), no. 2, 403-443.}
%
\bibitem{GG2}{
G. Gentile, G. Gallavotti,
\textit{Degenerate elliptic resonances},
Comm. Math. Phys.
\textbf{257} (2005), no. 2, 319--362.}
%
\bibitem{GM}{
G. Gentile, V. Mastropietro,
\textit{Methods for the analysis of the Lindstedt series
for KAM tori and renormalizability in classical mechanics.
A review with some applications},
Rev. Math. Phys.
\textbf{8} (1996), no. 3, 393--444. }
%
\bibitem{KL}{
H. Koch, J. Lopes Dias,
\textit{Renormalization of Diophantine skew flows,
with applications to the reducibility problem},
Preprint.}
%
\bibitem{K1}{
R. Krikorian,
\textit{R\'eductibilit\'e des syst\`emes produits-crois\'es \`a
valeurs dans des groupes compacts},
Ast\'erisque
\textbf{259} (1999), vi+216 pp.}
%
\bibitem{K2}{
R. Krikorian,
\textit{Global density of reducible quasi-periodic cocycles
on ${\bf T}^{1} \times {\rm SU}(2)$},
Ann. of Math.
\textbf{154} (2001), no. 2, 269--326.}
%
\bibitem{I}{
A. Iserles,
\textit{Expansions that grow on trees},
Notices Amer. Math. Soc.
\textbf{49} (2002), no. 4, 430--440.}
%
\bibitem{IN}{
A. Iserles, S.P. N{\o}rsett,
\textit{On the solution of linear differential equations in Lie groups},
R. Soc. Lond. Philos. Trans. Ser. A Math. Phys. Eng. Sci.
\textbf{357} (1999), no. 1754, 983--1019.}
%
\bibitem{JS}{
\`A. Jorba, C. Sim\'o,
\textit{On the reducibility of linear differential equations
with quasiperiodic coefficients},
J. Differential Equations
\textbf{98} (1992), no. 1, 111--124.}
%
\bibitem{L}{
J. Lopes Dias,
\textit{A normal form theorem for Brjuno skew-systems
through renormalization},
Preprint.}
%
\bibitem{MP}{
J. Moser, J. P\"oschel,
\textit{An extension of a result by Dinaburg and Sinai
on quasiperiodic potentials},
Comment. Math. Helv.
\textbf{59} (1984), no. 1, 39--85.}
%
\bibitem{PF}{
L. Pastur, A. Figotin,
\textit{Spectra of random and almost-periodic operators},
Grundlehren der Mathematischen Wissenschaften 297,
Springer, Berlin, 1992.}
%
\bibitem{R}{
H. R\"ussmann,
\textit{On the one-dimensional Schr\"odinger equation
with a quasiperiodic potential},
Nonlinear dynamics (International Conference, New York, 1979),
pp. 90--107, Ann. New York Acad. Sci.,
357, New York Acad. Sci., New York, 1980.}
%
\bibitem{S}{
W.M. Schmidt,
\textit{Diophantine approximation},
Lecture Notes in Mathematics, 785, Springer, Berlin, 1980.}
%
\bibitem{XZ}{
Ju. Xu, Q. Zheng,
\textit{On the reducibility of linear differential equations
with quasiperiodic coefficients which are degenerate},
Proc. Amer. Math. Soc.
\textbf{126}(1998), no. 5, 1445--1451.}

\end{thebibliography}
\end{document}